\documentstyle[chords,amssymb,epsf,12pt]{article}

\newtheorem{theor}{Theorem}[section]
\newtheorem{thm}[theor]{Theorem}
\newtheorem{lemma}[theor]{Lemma}
\newtheorem{prop}[theor]{Proposition}
\newtheorem{cor}[theor]{Corollary}

\newtheorem{dfn}[theor]{Definition}
\newtheorem{rem}[theor]{Remark}

\newcommand{\begproof}{\noindent \\ {\em Proof\/}. }
\newcommand{\eproof}{{\ \hfill  $\Box$}}

\newcommand{\be}{\begin{equation}}
\newcommand{\ee}{\end{equation}}
\newcommand{\bea}{\begin{eqnarray}}
\newcommand{\eea}{\end{eqnarray}}
\newcommand{\bean}{\begin{eqnarray*}}
\newcommand{\eean}{\end{eqnarray*}}
\newcommand{\abb}{\addtolength{\belowdisplayskip}{\belowdisplayskip}}
\newcommand{\half}{\textstyle {1\over 2}}
\newcommand{\text}[1]{ {\rm {#1} } }           
\renewcommand{\Bbb}{\mathbb}              

\newcommand{\Q}{{\Bbb Q}}
\newcommand{\C}{{\Bbb C}}
\newcommand{\R}{{\Bbb R}}
\newcommand{\Z}{{\Bbb Z}}

\renewcommand{\S}{\Sigma}
\newcommand{\cA}{{\cal A}}
\newcommand{\cB}{{\cal B}}
\newcommand{\cC}{{\cal C}}
\newcommand{\cF}{{\cal F}}
\newcommand{\cK}{{\cal K}}
\newcommand{\cD}{{\cal D}}
\newcommand{\cG}{{\cal G}}

\newcommand{\cU}{{\cal U}}

\newcommand{\cW}{{\cal W}}

\newcommand{\Ker}{\mathop{\rm Ker}}
\newcommand{\Mor}{\mathop{\rm Mor}}
\newcommand{\Inv}{\mathop{\rm Inv}}
\newcommand{\I}{\mathord{\rm I}}
\newcommand{\II}{\mathord{\rm I\!I}}
\newcommand{\III}{\mathord{\rm I\!I\!I}}
\newcommand{\str}{\mathop{\rm str}}
\newcommand{\tr}{\mathop{\rm tr}}
\newcommand{\BN}{Bar-Natan}
\newcommand{\Fn}{Feynman\ }
\newcommand{\Fd}{\Fn diagram}

\newcommand{\gl}{\mbox{$gl(1|1)$}}
\newcommand{\ra}{\rightarrow}

\newcommand{\Vas}{Vassiliev}
\newcommand{\VI}{\Vas\ invariants}
\newcommand{\ki}{knot invariant}

\newcommand{\ws}{weight system}

\newcommand{\ord}{{\mathop{\rm ord}}}

\begin{document}
\title{
\bf The universal Vassiliev invariant for the Lie superalgebra
$gl(1|1)$
}

\author{
Jos\'{e}~M~Figueroa-O'Farrill\thanks{e-mail: {\tt
j.m.figueroa@qmw.ac.uk}.\qquad
Supported by the EPSRC under contract GR/K57824.}
\\
{\em Department of Physics, Queen Mary and Westfield College},\\
{\em London E1 4NS, UK,}
 \\ \ \\
Takashi Kimura\thanks{e-mail: {\tt kimura@math.bu.edu}, \qquad This research
was supported in part by an NSF postdoctoral research fellowship}
\\ 
{\em Department of Mathematics, Boston University,}\\ 
{\em 111 Cummington Street; Boston, MA 02215, USA,}
\medskip\\ and \medskip\\
Arkady Vaintrob\thanks{e-mail: {\tt vaintrob@math.utah.edu}, \qquad
on leave from New Mexico State University, Las Cruces, NM 88003}
\\ {\em Department of Mathematics, University of Utah,} 
\\ {\em Salt Lake City, UT 84112, USA}
}

\date{January 1996}
\maketitle

\begin{abstract}
We compute the universal \ws\ for \VI\ coming from the Lie
superalgebra \gl\ applying the construction of \cite{YB}. This \ws\ 
is
a function from the space of chord diagrams to the center $Z$ of the
universal enveloping algebra of \gl, and we find a combinatorial
expression for it in terms of the standard generators of $Z$. 
The resulting \ki s generalize the Alexander-Conway polynomial.    
\end{abstract}
\newpage

\section*{Introduction}

Vassiliev in \cite{V1} initiated a study of a new class of knot
invariants which attracted a lot of interest because of their ability
to distinguish knots as well as their aesthetic beauty and numerous
connections with other classical (as well as {\em quantum}) fields of
mathematics and physics.

The space of Vassiliev invariants has a natural filtration $V_0\subset
V_1 \subset V_2 \ \cdots$ whose (adjoint) quotients $V_n/V_{n-1}$ can
be described in terms of combinatorial objects called \ws s (as it was
shown in the works of Birman--Lin \cite{BL}, Bar-Natan \cite{BN} and
Kontsevich \cite{K}).  Weight systems of order $n$ are functions on
the set $\cD_n$ of chord diagrams --- circles with $n$ chords
(unordered pairs of points) satisfying certain relations.  They form
a finite-dimensional vector space $W_n$ and to describe this space
(and thus to find all Vassiliev invariants of order $n$) we just need
to solve a system of linear equations.  But the numbers of unknowns
and equations in this system grow extremely fast with $n$
(cf. \cite{BN}) and therefore we need a better way to approach $W_n$.

Motivated by perturbative Chern--Simons theory, Bar-Natan
\cite{BNCS} and Kon\-tse\-vich \cite{K} gave a construction of weight
systems using a Lie algebra $L$ with an invariant inner product and a
module.  In \cite{YB} one of us generalized this construction to
so-called self-dual Yang--Baxter Lie algebras.

This class of algebras includes Lie algebras and Lie superalgebras
with invariant inner products.  Each such algebra produces a sequence
of \ws s $W_L: \cD_n \rightarrow \cU(L)$ with values in the center
$Z_L$ of the universal enveloping algebra $\cU(L)$.  We call this
family of \ws s a {\em universal weight system\/} for $L$, since \ws s
corresponding to $L$-modules are obtained from $W_L$ by taking traces
in the corresponding representation.

In \cite{BNCS,BN} Bar-Natan computed \ws s for the defining
representations of classical Lie algebras.  Although for $L=sl_n$ the
center of $\cU(L)$ is dual to the space of functionals on $\cU$
spanned by the characters of the exterior powers of the defining
representations and the corresponding weight systems are known
(cf. \cite{BN}), the problem of finding a direct combinatorial
expression for $W_L$ in terms of standard generators of $Z_L$ is
highly nontrivial.

The first universal weight system was studied by Chmutov and Varchenko
\cite{CV} who considered the case where $L = sl_2$.  The center of
$\cU(sl_2)$ is isomorphic to the polynomial algebra $\C[c]$, where $c$ is the
Casimir element of $\cU(sl_2)$, and the main result of \cite{CV} is a
recursive formula for computing values of $W_{sl_2}$ on chord diagrams.

In this paper we are studying the universal weight system for \gl, the
simplest interesting example of a self-dual Lie superalgebra.

The case of \gl\ is different from that of $sl_2$ in several ways.
First, the center of $\cU(gl(1|1))$ has two generators $h$ and $c$.
Second, all the invariant functionals on $\cU(sl_2)$ are linear
combinations of traces in irreducible representations, which is not
true for \gl\ (since its superdimension as well as the superdimension
of its generic modules is $0$).  Third, a choice of values for $h$ and
$c$ in our weight system gives rise to the sequence of coefficients of
the Alexander-Conway polynomial.  This once again confirms the
relationship between \gl\ supersymmetry and the Alexander polynomial
(cf. \cite{KS,RS1}). Weight systems coming from Lie superalgebras
other than \gl\ can be used to explain relations between classical
knot invariants (cf. \cite{YB,algstr}) and to construct invariants
that are stronger than all invariants corresponding to semi-simple Lie
algebras (cf. \cite{Vogel}).
\medskip

Let us formulate our answer.  Let $c$ be the quadratic Casimir for
\gl, and $h\in \gl$, the identity matrix. Then $Z(\cU(\gl))= 
\C[h,c]$, and our main result is the following recursive formula for
values of the weight system $W=W_{gl(1|1)}$.
\medskip

\noindent
{\bf Theorem.} {\em
Let $D$ be a chord diagram, ``$a$'' a fixed chord in $D$ and
$b_1,b_2,\ldots$ are all the chords of $D$ intersecting $a$. Denote by
$D_a$ (resp.  $D_{a,i}$, $D_{a,ij}$) the diagram $D - a$ (resp. $D - a
- b_i$, $D- a - b_i - b_j$). Then
$W(D)$ is a polynomial in $c$ and $y=-h^2$ satisfying
\begin{eqnarray}
W(D) &=& c \, W(D_a) - y  \sum_i W(D_{a,i})    \label{eq:rec} \\
&&{}+ y \sum_{i<j} \left(W\left(D_{a,ij}^{+-}\right) +
W\left(D_{a,ij}^{-+}\right)- W\left(D_{a,ij}^l\right) -
W\left(D_{a,ij}^r\right)\right)~, \nonumber
\end{eqnarray}
where $D_{a,ij}^{+-}$ (resp. $D_{a,ij}^{-+}$; $D_{a,ij}^l$;
$D_{a,ij}^r$) is the diagram obtained by adding to $D_{a,ij}$ a new
chord connecting the left end of $b_i$ and the right end of $b_j$
(resp. the right end of $b_i$ and the left end of $b_j$; the left ends
of $b_i$ and $b_j$; the right ends of $b_i$ and $b_j$) assuming that
the chord $a$ is drawn vertically. Pictorially, if $D$ denotes the
diagram below and $a$ is the vertical chord and $i,j$ are chords which
intersect $a$ so that $i$ is the upper chord and $j$ is the lower
chord then
\begin{equation}
\begin{array}{rrrr}
D=\Picture{
\DottedCircle\FullChord[1,5]\FullChord[3,9]\FullChord[7,11]
\put(0.18,0.05){\makebox(0,0){$^a$}}
\put(-0.35,0.55){\makebox(0,0){$^i$}}
\put(-0.35,-0.25){\makebox(0,0){$^j$}}
}&
D_a= \Picture{
\DottedCircle\FullChord[1,5] \FullChord[7,11]\Arc[3]\Arc[9]
}&
D_{a,i} = \Picture{
\DottedCircle\FullChord[7,11]\Arc[3]\Arc[9]\Arc[1]\Arc[5]
}&
D_{a,j}=\Picture{
\DottedCircle\FullChord[1,5]\Arc[3]\Arc[9]\Arc[7]\Arc[11]
} \\ [30pt]
D_{a,ij}^{+-} = \Picture{
\DottedCircle\FullChord[5,11]\Arc[1] \Arc[3]\Arc[7]\Arc[9]
}&
D_{a,ij}^{-+} = \Picture{
\DottedCircle\FullChord[1,7]\Arc[3]\Arc[5]\Arc[9]\Arc[11]
}&
D_{a,ij}^{l} = \Picture{
\DottedCircle\FullChord[5,7]\Arc[1]\Arc[3]\Arc[9]\Arc[11]
}&
D_{a,ij}^{r} = \Picture{
\DottedCircle\FullChord[1,11]\Arc[3] \Arc[5]\Arc[7]\Arc[9]
}
\end{array}
\label{eq:pixD}
\end{equation}
}
\medskip

Since all the diagrams in the right-hand side of (\ref{eq:rec}) have
less chords than $D$, this allows us to compute the value of the
weight system $W_{gl(1|1)}$ on any chord diagram recursively (and
quite effectively); see the examples in subsection \ref{sec:recursion}
and the tables in the appendix.

\bigskip

The outline of the paper is as follows.  In Section 1 we collect
preliminary facts about \VI, chord diagrams and \ws s, and recall the
construction from \cite{YB} of \ws s based on Lie superalgebras (or,
more generally, on Yang--Baxter Lie algebras).

In Section 2 we present the necessary algebraic information about \gl,
its representations, invariant tensors and the center of the universal
enveloping algebra.  An identity between invariant tensors of fourth
order, which is pivotal for proving formula (\ref{eq:rec}), is also
established here.

In Section 3 we prove the recursive formula (\ref{eq:rec}) for the
universal weight system $W_{gl(1|1)}$.
We also prove that deframing the universal weight system consists
simply in evaluating it at $c{=}0$.
We conclude the section with an alternative way of producing the same
weight system.

In Section 4 we show how the Alexander-Conway polynomial can be
obtained from our \ws.

\bigskip

\noindent
{\bf Acknowledgments.}  
The work described in this paper has been done when the first and the
third authors visited the second at the University of North Carolina
at Chapel Hill.  We are all grateful to the Department of Mathematics
of UNC and to Jim Stasheff for their hospitality.  Additional thanks
for hospitality go from JMF to Louise Dolan and the Department of
Physics of UNC and from AV to the Max-Planck Institut f\"ur
Mathematik where a large part of the paper was written up.
It is a pleasure to thank Bill Spence for helpful conversations.

We would like to convey our gratitude to Craig Snydal without whom
this paper would not have been written. In his thesis \cite{Craig},
Craig wrote computer programs which computed universal Vassiliev
invariants associated to $gl(1|1)$ up to order $5$ before the
recursion relation was discovered. It was while supervising Craig's
thesis that TK learned a great deal about this subject. In particular,
the fundamental relation (\ref{eq:tensoreq}) was initially discovered
using software written by TK which was incorporated into Craig's
programs. Craig's computations proved to be useful in checking our
results, as well. To him we give our special thanks.


\section{Vassiliev invariants and Yang--Baxter Lie algebras}

Here we recall some concepts and results related to Vassiliev
invariants which we will need later. In particular, we 
review their relationship with Lie algebra-type structures. For more
details see \cite{BN,K,YB}.
 
\subsection{Vassiliev invariants, chord diagrams and weight systems}

\begin{dfn}{\rm
A {\em singular knot\/} is an immersion $K: S^1 \rightarrow \R^3$ with
a finite number of transversal double self-intersections (or double
points).  The set of singular knots with $n$ double points is denoted
by $\cK_n$.

A {\em chord diagram\/} of order $n$ is an oriented circle with $n$
non-intersecting pairs of points (chords) on it, up to an orientation
preserving diffeomorphism of the circle.  Denote by $\cD_n$ the set of
all chord diagrams with $n$ chords.

Every $K\in \cK_n$ has a chord diagram $ch(K)\in \cD_n$ whose chords
are the inverse images of the double points of $K$.  } \end{dfn}

Vassiliev showed that every  knot invariant $I$ with values in an abelian
group $k$ extends to an invariant of singular knots by the rule
\be 
I(K_0) = I(K_+) - I(K_-),    \label{eq:vasrel}
\ee
where $K_0$, $K_+$, and $K_-$ are singular knots which differ only
inside a small ball as shown below:

\def\KZero{
\begin{picture}(2,2)(-1,-1)
\put(0,0){\circle{2}}
\put(-0.707,-0.707){\vector(1,1){1.414}}
\put(0.707,-0.707){\vector(-1,1){1.414}}
\put(0,0){\circle*{0.15}}
\end{picture}}

\def\KPlus{
\begin{picture}(2,2)(-1,-1)
\put(0,0){\circle{2}}
\put(-0.707,-0.707){\vector(1,1){1.414}}
\put(0.707,-0.707){\line(-1,1){0.6}}
\put(-0.107,0.107){\vector(-1,1){0.6}}
\end{picture}}

\def\KMinus{
\begin{picture}(2,2)(-1,-1)
\put(0,0){\circle{2}}
\put(-0.707,-0.707){\line(1,1){0.6}}
\put(0.107,0.107){\vector(1,1){0.6}}
\put(0.707,-0.707){\vector(-1,1){1.414}}
\end{picture}}

\begin{equation}
\mathop{\KZero}_{K_0}\qquad
\mathop{\KPlus}_{K_+}\qquad
\mathop{\KMinus}_{K_-}\label{eq:VasRel}
\end{equation}

\begin{dfn}{\rm
A knot invariant $I$ is called an {\em invariant of order\/} ($\le$)
$n$ if $I(K)=0$ for any $K\in \cK_{n+1}$.

The $k$-module of all invariants of order $n$ is denoted by $V_n$.  We
have the obvious filtration
\begin{displaymath}
V_0 \subset V_1 \subset V_2 \cdots
\subset V_n \subset \cdots .
\end{displaymath}

Elements of $V_\infty = \cup_n V_n$ are called {\em invariants of
finite type\/} or {\em Vassiliev invariants}.

Analogously, Vassiliev invariants of framed knots (and links) can be
defined.  }\end{dfn}

\medskip

An immediate corollary of the definition of Vassiliev invariants is
that the value of an invariant $I \in V_n$ on a singular knot $K$ with
$n$ self-intersections depends only on the diagram $ch(K)$ of $K$.  In
other words, $I$ descends to a function on $\cD_n$ which we will still
denote by $I$ slightly abusing notation.  These functions satisfy two
groups of relations.
\begin{equation}
I\left(\Picture{\DottedCircle\FullChord[2,10]\Arc[1]\Arc[0]\Arc[11]}\right)
= 0\label{eq:1term}
\end{equation}
and
\begin{equation}
I\left(
\Picture{\DottedCircle\FullChord[1,8]\Arc[2]\FullChord[5,9]}
\right) -
I\left(
\Picture{\DottedCircle\FullChord[1,9]\Arc[2]\FullChord[5,8]}
\right) +
I\left(
\Picture{\DottedCircle\FullChord[2,5]\Arc[8]\FullChord[1,9]}
\right) -
I\left(
\Picture{\DottedCircle\FullChord[1,5]\Arc[8]\FullChord[2,9]}
\right) = 0\label{eq:4term}
\end{equation}

(When drawing a chord diagram we always assume that the circle is
oriented counterclockwise.  Chords whose endpoints lie on the solid
arcs are shown explicitly.  The diagram may contain other chords whose
endpoints lie on the dotted arcs, provided that they are the {\em
same\/} in all diagrams appearing in the same relation.)

\begin{dfn}{\rm
A function $W: {\cal D}_n \rightarrow k$ is called a {\em weight
system of order $n$\/} if it satisfies the conditions
(\ref{eq:4term}). If, in addition, $W$ satisfies the relations
(\ref{eq:1term}), then we call it a {\em strong \ws}.  } \end{dfn}

Denote by $\cW_n$ (resp. by $\overline{\cW}_n$) the set of all weight
systems (resp. strong weight systems) of order $n$.

It is easy to see that the natural map $V_n/V_{n-1} \rightarrow
\overline{\cW}_n$ is injective.

The remarkable fact proved by Kontsevich and Bar-Natan is that this
map is also surjective (at least when $k \supset \Q$).  In other
words, each strong weight system of order $n$ is a restriction to
${\cal{D}}_n$ of some Vassiliev invariant.

Let ${\cal{A}}_n$ (resp.  $\bar{\cal{A}}_n$ ) be the dual space to
$\cW_n$ (resp. $\overline{\cW}_n$), i.e.,  the space of formal linear
combinations of diagrams from ${\cD_n}$ modulo four-term
relations\footnote{which are now considered as relations in the space
of diagrams, i.e., without $I$.}{} (\ref{eq:4term}) (resp. four- and
one-term (\ref{eq:1term}) relations), and
\begin{displaymath}
{\cal{A}}=\oplus_n {\cal{A}}_n,
\quad
\bar{\cal{A}}=\oplus_n\bar{\cal{A}}_n.
\end{displaymath}

Denote by $V^f_n$ the space of Vassiliev invariants of order $\le n$
of {\em framed\/} knots.

\begin{thm}{\cite{BN,K}.} \label{thm:cdalg}
If $k\supset \Q$ 
then

(1) $V_n/V_{n-1} \simeq \overline{\cW}_n \simeq \bar{\cA}_n^*.$

\smallskip

(2) $V^f_n/V^f_{n-1} \simeq \cW_n \simeq \cA_n^*.$

(3) The operation of connected sum of diagrams induces on $\cal{A}$
and $\bar{\cal{A}}$ structures of commutative graded $k$-algebras.
The comultiplication dual to this product makes the graded algebras
$\cW =\oplus_n \cW_n \simeq {\cal{A}}^*$ and $\overline{\cW}=\oplus_n
\overline{\cW}_n \simeq \bar{\cal{A}}^*$ into commutative and
co-commutative Hopf algebras.

(4) ${\cal{A}} = \bar{\cal{A}}[\Theta]$, where $\Theta$ is the
primitive element in ${A}_1$ represented by the only chord diagram of
order $1$, i.e.,
\begin{displaymath}
\cA_n=\bar{\cA_0}\otimes\Theta^n+\bar{\cA_1}\otimes\Theta^{n-1}
+\ldots \bar{\cA_n}\otimes 1~.
\end{displaymath}
\end{thm}

Part (4) of the theorem shows that there is a canonical projection
(``deframing'') $\cW_n \rightarrow \overline{\cW}_n$.  Therefore, for
every weight system there is a canonical strong weight system (and
consequently, a knot invariant).
We will not be concerned with the one-term relations in the bulk of
the paper, but will return to it briefly in Section 3.4 when we
discuss the deframed universal weight system.

\begin{rem}
{\rm
Similar results hold for links. The only difference is that in
this case we will have to consider several circles.
}\end{rem}

\subsection{Weight systems coming from Yang--Baxter algebras  } 
\label{sec:FG}

Here we recall a construction from \cite{YB} that assigns a family of
weight systems to every Yang--Baxter Lie algebra with an invariant
inner product.

First, we introduce a more general class of diagrams.

\begin{dfn}{\rm
A {\em \Fd \/} of order $p$ is a graph with $2p$ vertices of degrees 1
or 3 with a cyclic ordering on the set of its univalent ({\em
external\/}) vertices and on each set of $3$ edges meeting at a
trivalent ({\em internal\/}) vertex.\footnote{A Feynman diagram after
forgetting the ordering of external vertices is called a Chinese
Character diagram in \cite{BN}.}  Let $\cF_p$ denote the set of all
\Fd s with $2p$ vertices (up to the natural equivalence of graphs with
orientations). The set $\cD_p$ of chord diagrams with $p$ chords is a
subset of $\cF_p$.

\medskip

We draw \Fd s by placing their external vertices ({\em legs}) on a
circle which is oriented in the counterclockwise direction called the
{\em Wilson loop\/}.  The edges of a \Fd\ are called {\em
propagators}.  We assume that the propagators meeting at each internal
vertex are oriented counterclockwise.  }
\end{dfn}

Denote by ${\cal{B}}_p$ the vector space generated by \Fd s of order
$p$ modulo relations

\begin{equation}\label{eq:3term}
\Picture{
\Arc[9]
\Arc[10]
\Arc[11]
\Arc[0]
\DottedArc[1]
\DottedArc[8]
\thinlines
\put(0.7,-0.7){\circle*{0.15}}
\qbezier(-0.1,0.1)(0.2,-0.5)(0.5,0.1)
\qbezier(0.2,-0.2)(0.3,-0.6)(0.7,-0.7)
\put(0,-1.4){\makebox(0,0){${}_{D_Y}$}}
}\quad = \quad
\Picture{
\Arc[9]
\Arc[10]
\Arc[11]
\Arc[0]
\DottedArc[8]
\DottedArc[1]
\Endpoint[10]
\Endpoint[11]
\thinlines
\put(0.5,-0.866){\line(-3,5){0.58}}
\put(0.866,-0.5){\line(-3,5){0.36}}
\put(0,-1.4){\makebox(0,0){${}_{D_{\mid\mid}}$}}
} \quad - \quad 
\Picture{
\Arc[9]
\Arc[10]
\Arc[11]
\Arc[0]
\DottedArc[8]
\DottedArc[1]
\Endpoint[10]
\Endpoint[11]
\thinlines
\put(0.5,0.1){\line(0,-1){0.97}}
\put(0.866,-0.5){\line(-5,3){1}}
\put(0,-1.4){\makebox(0,0){${}_{D_X}$}}
}
\end{equation}

\ 

\bigskip

More precisely, $\cB_p = \langle\cF_p\rangle /\langle
D_Y-D_{\mid\mid}-D_X\rangle $, where the diagrams $D_{\mid\mid}$ and
$D_X$ are obtained from the diagram $D_Y$ by replacing its
$Y$-fragment by the $\mid\mid$- and $X$- fragments respectively.
This gives another description of the space
$(V^f_p/V^f_{p-1})^*=\cA_p$.

\begin{prop}{\cite{BN}}  \label{prop:fd}

(1) The embedding $\cD_p \hookrightarrow \cF_p$ induces an isomorphism
\quad $\cB_p \simeq \cA_p.$

(2) The following local relations hold for internal vertices in \Fd s:

\begin{displaymath}
\begin{picture}(2,2)(0,0.375)
\qbezier(0.5,2)(1.65,1.3)(1.5,1)
\qbezier(1.5,2)(0.35,1.3)(0.5,1)
\qbezier(0.5,1)(1,0)(1.5,1)
\put(1,0){\line(0,1){0.5}}
\end{picture}
= - 
\begin{picture}(2,2)(0,0.375)
\qbezier(0.5,2)(1,0)(1.5,2)
\put(1,0){\line(0,1){1}}
\end{picture}
\qquad{\mathrm and}\qquad
\begin{picture}(3,2)(-1,-1)
\put(0,-1.4){\line(1,1){1.4}}
\qbezier(-0.1,0.1)(0.2,-0.5)(0.5,0.1)
\qbezier(0.2,-0.2)(0.3,-0.6)(0.7,-0.7)
\end{picture} = 
\begin{picture}(3,2)(-1,-1)
\put(-0.034,-1.4){\line(1,1){1.4}}
\put(0.5,-0.866){\line(-3,5){0.58}}
\put(0.866,-0.5){\line(-3,5){0.36}}
\end{picture} - 
\begin{picture}(3,2)(-1,-1)
\put(-0.034,-1.4){\line(1,1){1.4}}
\put(0.5,0.1){\line(0,-1){0.97}}
\put(0.866,-0.5){\line(-5,3){1}}
\end{picture}
\end{displaymath}
\end{prop}
\smallskip

\begin{dfn}{\rm
A {\em Feynman graph\/} is a graph with $1$- and $3$-valent vertices
such that the set of its legs (univalent vertices) is a disjoint union
of two {\em linearly\/} ordered sets: {\em incoming\/} and {\em
outgoing\/} legs and the set of propagators meeting at each of its trivalent
({\em internal\/}) vertices has a cyclic ordering.

We denote by ${\cF}_{a,b}$ the set of Feynman graphs with $a$ incoming
and $b$ outgoing legs and by $\cF_*$ the set of all Feynman graphs.  }
\end{dfn}

We draw Feynman graphs by putting all incoming legs on a horizontal
segment, outgoing legs on a parallel segment below, and all the
internal vertices between these two segments (with the
counterclockwise orientation at each vertex).  Notice that closing
these segments into a circle, we turn a Feynman graph into a \Fd.
\medskip

Feynman graphs can be regarded as morphisms in a tensor category
$\bf FG$ with objects $0,1,2,\ldots$ and natural operations of
composition and tensor product.

If $A\in \cF_{b,c}$ and $B\in \cF_{a,b}$, then their composition
$A\circ B \in \cF_{a,c}$ is the Feynman graph obtained by attaching
the outgoing vertices of $B$ to the corresponding incoming vertices of
$A$.

The tensor product of $A\in \cF_{a,b}$ and $C\in \cF_{c,d}$ is the
graph $A\otimes C \in \cF_{a+c,b+d}$ obtained by placing $C$ to the
right of $A$.
\medskip

A natural way to assign invariants to Feynman graphs (and, ultimately,
to Feynman and chord diagrams) is to consider representations of the
category ${\bf FG}$. Recall that a {\em representation\/} of a tensor
category $\cC$ is a tensor functor $F$ from $\cC$
to the category of vector
spaces, i.e., an assignment of a vector space $F(A)$ to each object $A$
of $\cC$ and a linear map $F(g): F(A) \rightarrow F(B)$ to each
morphism $g\in \Mor(A,B)$ such that it respects composition and tensor
products.

In the case $\cC = {\bf FG}$ a representation is specified by a choice
of a vector space $L$ and a $k$-linear map $F(\Gamma) : L^{\otimes a}
\rightarrow L^{\otimes b}$ for every $\Gamma \in \cF_{a,b}$ (where
$L^{\otimes 0} = k$).

Every element of $\cF_*=\Mor({\bf FG})$ can be obtained by means of
composition and tensor product from the following {\em elementary
graphs\/}:
\bigskip
\begin{displaymath}
\begin{array}{c@{\hspace{30pt}}c@{\hspace{30pt}}c}
I=
\begin{picture}(2,2)(0,1)
\put(1,0){\line(0,1){2}}
\thicklines
\put(0,0){\line(1,0){2}}
\put(0,2){\line(1,0){2}}
\end{picture}
\in \cF_{1,1}, & 
b=
\begin{picture}(2,2)(0,1)
\qbezier(0.5,2)(1,0)(1.5,2)
\thicklines
\put(0,2){\line(1,0){2}}
\put(0,0){\line(1,0){2}}
\end{picture}
\in \cF_{2,0}, &
c= 
\begin{picture}(2,2)(0,1)
\qbezier(0.5,0)(1,2)(1.5,0)
\thicklines
\put(0,0){\line(1,0){2}}
\put(0,2){\line(1,0){2}}
\end{picture}
\in \cF_{0,2}, \\ [20pt]
f=
\begin{picture}(2,2)(0,1)
\qbezier(0.5,2)(1,0)(1.5,2)
\put(1,0){\line(0,1){1}}
\thicklines
\put(0,0){\line(1,0){2}}
\put(0,2){\line(1,0){2}}
\end{picture}
\in \cF_{2,1}, &
g=
\begin{picture}(2,2)(0,1)
\qbezier(0.5,0)(1,2)(1.5,0)
\put(1,2){\line(0,-1){1}}
\thicklines
\put(0,0){\line(1,0){2}}
\put(0,2){\line(1,0){2}}
\end{picture}
\in \cF_{1,2}, & 
S=
\begin{picture}(2,2)(0,1)
\put(0.5,0){\line(1,2){1}}
\put(1.5,0){\line(-1,2){1}}
\thicklines
\put(0,0){\line(1,0){2}}
\put(0,2){\line(1,0){2}}
\end{picture}
\in \cF_{2,2}.
\end{array}
\end{displaymath}
\
\bigskip

\noindent
Therefore, to specify a representation of $\bf FG$ we need to fix a
vector space $L$ and six tensors corresponding to the generators $I,
b, c, f, g,$ and $S$.

In \cite{YB} the defining relations between these generators were
found and it was shown that representations of the tensor category
$\bf FG$ on the category of vector spaces are in one-to-one
correspondence with the set of self-dual Yang--Baxter algebras.
\medskip

\begin{dfn}{\rm
A vector space $L$ over a field $k$ with an operator 
$$S: L\otimes L \rightarrow L\otimes L$$ 
and a multiplication $f: L\otimes L\rightarrow
L$ satisfying conditions (1-3) below is called a {\em Yang--Baxter
algebra\/} or, simply, an $S$-algebra. If, in addition, $L$ is
finite-dimensional and is equipped with an inner product $b: L\otimes
L \rightarrow k$ satisfying (4-7) below, it is called a {\em symmetric
self-dual $S$-algebra\/}\footnote{This nomenclature is not standard.
Such algebras have been termed ``Euclidean'' in \cite{YB} and
``metric'' in \cite{BN,BNCS}.  Throughout this paper we will refer to
such algebras simply as self-dual, leaving implicit the fact that $b$
is $S$-symmetric.}.  A {\em (self-dual) Yang--Baxter Lie algebra\/} is
a (self-dual) $S$-algebra $L$, such that the multiplication $f:
L\otimes L\rightarrow L$ obeys (8-9) below.

We use the standard notation: if $T\in {\text{Hom}}(V\otimes V, W)$,
then $T_{12}$ denotes the operator 
$$T\otimes {\text{id}}: V\otimes V \otimes U \rightarrow W\otimes U,$$
 etc.

\begin{center}

\begin{tabular}{@{}l@{.\hspace{6pt}}p{180pt}@{\hspace{10pt}}p{180pt}@{}}
1&\raggedright The operator $S$ is a {\em symmetry\/}:
\[S^2=id_{L\otimes L}\]
\vspace{-20pt}&
\vspace{-26pt}
\[
\begin{picture}(2,2)(0,0.775)
\qbezier(0.5,0)(2,1)(0.5,2)
\qbezier(1.5,0)(0,1)(1.5,2)
\thicklines
\put(0,0){\line(1,0){2}}
\put(0,2){\line(1,0){2}}
\end{picture}
=
\begin{picture}(2,2)(0,0.775)
\put(0.5,0){\line(0,1){2}}
\put(1.5,0){\line(0,1){2}}
\thicklines
\put(0,0){\line(1,0){2}}
\put(0,2){\line(1,0){2}}
\end{picture}
\]
\\ [9pt]
2&\raggedright $S$ satisfies the quantum Yang--Baxter equation:
\[S_{12}S_{23}S_{12}=S_{23}S_{12}S_{23}\]
\vspace{-20pt}&
\vspace{-12pt}
\[
\begin{picture}(3,2)(0,0.775)
\qbezier(1.5,0)(0,1)(1.5,2)
\put(0.5,0){\line(1,1){2}}
\put(0.5,2){\line(1,-1){2}}
\thicklines
\put(0,0){\line(1,0){3}}
\put(0,2){\line(1,0){3}}
\end{picture}
=
\begin{picture}(3,2)(0,0.775)
\qbezier(1.5,0)(3,1)(1.5,2)
\put(0.5,0){\line(1,1){2}}
\put(0.5,2){\line(1,-1){2}}
\thicklines
\put(0,0){\line(1,0){3}}
\put(0,2){\line(1,0){3}}
\end{picture}
\]
\end{tabular}
\end{center}

\begin{center}
\begin{tabular}{@{}l@{.\hspace{6pt}}p{180pt}@{\hspace{10pt}}p{180pt}@{}}
3&\raggedright The multiplication $f$ is compatible with the symmetry
$S$:
\[S\,f=f_{23}S_{12}S_{23}\]
\vspace{-20pt}&
\vspace{-12pt}
\[
\begin{picture}(3,2)(0,0.775)
\qbezier(1,2)(1.5,0)(2,2)
\put(1.5,0){\line(0,1){1}}
\qbezier(2.5,2)(2.5,1)(0.5,0)
\thicklines
\put(0,0){\line(1,0){3}}
\put(0,2){\line(1,0){3}}
\end{picture}
=
\begin{picture}(3,2)(0,0.775)
\qbezier(1,2)(1.5,0)(2,2)
\put(1.5,0){\line(0,1){1}}
\qbezier(2.5,2)(0.5,1)(0.5,0)
\thicklines
\put(0,0){\line(1,0){3}}
\put(0,2){\line(1,0){3}}
\end{picture}
\]
\\ [9pt]
4&\raggedright The bilinear form $b$ is $S$-symmetric:
\[bS=S\]
\vspace{-20pt}&
\vspace{-12pt}
\[
\begin{picture}(2,2)(0,0.775)
\qbezier(0.5,2)(1.65,1.3)(1.5,1)
\qbezier(1.5,2)(0.35,1.3)(0.5,1)
\qbezier(0.5,1)(1,0)(1.5,1)
\thicklines
\put(0,0){\line(1,0){2}}
\put(0,2){\line(1,0){2}}
\end{picture}
=
\begin{picture}(2,2)(0,0.775)
\qbezier(0.5,2)(1,0)(1.5,2)
\thicklines
\put(0,0){\line(1,0){2}}
\put(0,2){\line(1,0){2}}
\end{picture}
\]
\end{tabular}
\end{center}

\begin{center}
\begin{tabular}{@{}l@{.\hspace{6pt}}p{180pt}@{\hspace{10pt}}p{180pt}@{}}
5&\raggedright $b$ is compatible with $S$:
\[b_{12}S_{23}=b_{23}S_{12}\]
\vspace{-20pt}&
\vspace{-26pt}
\[
\begin{picture}(3,2)(0,0.775)
\qbezier(0.5,2)(1.5,0)(2.5,2)
\qbezier(1.5,2)(0.5,1.5)(0.5,0)
\thicklines
\put(0,0){\line(1,0){3}}
\put(0,2){\line(1,0){3}}
\end{picture}
=
\begin{picture}(3,2)(0,0.775)
\qbezier(0.5,2)(1.5,0)(2.5,2)
\qbezier(1.5,2)(2.5,1.5)(2.5,0)
\thicklines
\put(0,0){\line(1,0){3}}
\put(0,2){\line(1,0){3}}
\end{picture}
\]
\end{tabular}
\end{center}

\begin{center}
\begin{tabular}{@{}l@{.\hspace{6pt}}p{180pt}@{\hspace{10pt}}p{180pt}@{}}
6&\raggedright $b$ is a non-degenerate bilinear form with inverse 
$c$:
\[b_{23}c_{12} = id_L = b_{12}c_{23}\]
\vspace{-20pt}&
\vspace{-12pt}
\[
\begin{picture}(2,2)(0,0.775)
\qbezier(0.5,2)(0.5,0)(1,1)
\qbezier(1.5,0)(1.5,2)(1,1)
\thicklines
\put(0,0){\line(1,0){2}}
\put(0,2){\line(1,0){2}}
\end{picture}
=
\begin{picture}(1,2)(0,0.775)
\put(0.5,0){\line(0,1){2}}
\thicklines
\put(0,0){\line(1,0){1}}
\put(0,2){\line(1,0){1}}
\end{picture}
=
\begin{picture}(2,2)(0,0.775)
\qbezier(0.5,0)(0.5,2)(1,1)
\qbezier(1.5,2)(1.5,0)(1,1)
\thicklines
\put(0,0){\line(1,0){2}}
\put(0,2){\line(1,0){2}}
\end{picture}
\]
\end{tabular}
\end{center}

\begin{center}
\begin{tabular}{@{}l@{.\hspace{6pt}}p{180pt}@{\hspace{10pt}}p{180pt}@{}}
7&\raggedright $b$ is $f$-invariant:
\[b\circ f_{12} = b\circ f_{23}\]
\vspace{-20pt}&
\vspace{-26pt}
\[
\begin{picture}(3,2)(0,0.775)
\qbezier(0.5,2)(1,1)(1.5,2)
\qbezier(1,1.5)(1.75,0)(2.5,2)
\thicklines
\put(0,0){\line(1,0){3}}
\put(0,2){\line(1,0){3}}
\end{picture}
=
\begin{picture}(3,2)(0,0.775)
\qbezier(1.5,2)(2,1)(2.5,2)
\qbezier(0.5,2)(1.25,0)(2,1.5)
\thicklines
\put(0,0){\line(1,0){3}}
\put(0,2){\line(1,0){3}}
\end{picture}
\]
\end{tabular}
\end{center}

\begin{center}
\begin{tabular}{@{}l@{.\hspace{6pt}}p{180pt}@{\hspace{10pt}}p{180pt}@{}}
8&\raggedright 
$f$ is $S$-{\em skew-symmetric\/}:
\[f\circ S=-f\]
\vspace{-20pt}&
\vspace{-26pt}
\[
\begin{picture}(2,2)(0,0.775)
\qbezier(0.5,2)(1.65,1.3)(1.5,1)
\qbezier(1.5,2)(0.35,1.3)(0.5,1)
\qbezier(0.5,1)(1,0)(1.5,1)
\put(1,0){\line(0,1){0.5}}
\thicklines
\put(0,0){\line(1,0){2}}
\put(0,2){\line(1,0){2}}
\end{picture}
= -
\begin{picture}(2,2)(0,0.775)
\qbezier(0.5,2)(1,0)(1.5,2)
\put(1,0){\line(0,1){1}}
\thicklines
\put(0,0){\line(1,0){2}}
\put(0,2){\line(1,0){2}}
\end{picture}
\]
\\ [10pt]
9&\raggedright $f$ satisfies the $S$-{\em Jacobi identity\/}:
\[f\circ f_{12} = f\circ f_{23} - f\circ f_{23}\circ S_{12}\]
\vspace{-20pt}&
\vspace{-26pt}
\[
\begin{picture}(2.5,2)(0,0.775)
\qbezier(0.5,2)(1,1)(1.5,2)
\qbezier(1,1.5)(1.5,0)(2,2)
\put(1.44,0){\line(0,1){0.865}}
\thicklines
\put(0,0){\line(1,0){2.5}}
\put(0,2){\line(1,0){2.5}}
\end{picture}
=
\begin{picture}(2.5,2)(0,0.775)
\qbezier(1,2)(1.5,1)(2,2)
\qbezier(0.5,2)(1,0)(1.5,1.5)
\put(1.06,0){\line(0,1){0.865}}
\thicklines
\put(0,0){\line(1,0){2.5}}
\put(0,2){\line(1,0){2.5}}
\end{picture}
-
\begin{picture}(2.5,2)(0,0.775)
\qbezier(0.5,2)(1.25,0)(2,2)
\qbezier(0.25,1)(0.75,0)(1.25,1)
\qbezier(0.25,1)(0.05,1.4)(1.5,2)
\put(0.75,0){\line(0,1){0.5}}
\thicklines
\put(0,0){\line(1,0){2.5}}
\put(0,2){\line(1,0){2.5}}
\end{picture}
\]
\\
\end{tabular}
\end{center}
} 
\end{dfn}

\medskip

Every self-dual $S$-algebra $L$ gives a representation of the category
${\bf FG}$ on the vector space $L$ (or, more accurately, a tensor
functor $F_L$ from $\bf FG$ to the category ${\cal{T}}(L)$ of tensor
powers of $L$).

Lie $S$-algebras arise when we attempt to construct representations of
${\bf FG}$ satisfying relations (2) of Proposition \ref{prop:fd}.
\medskip

Denote by $\cG_{a,b}$ the quotient space of the space $\langle
\cF_{a,b}\rangle $ of formal linear combinations of elements of
$\cF_{a,b}$ by relations \ref{prop:fd}(2) for internal vertices.

It is clear that $\cF_*$ induces on $\cG_* = \bigcup \cG_{a,b}$ the
structure of a tensor category with operations extended by
linearity.
\medskip

\begin{prop}
There is a one-to-one correspondence between representations of the
category $\cG_*$ and self-dual Yang--Baxter Lie algebras.
\end{prop}
\medskip

To go from the category $\cG_*$ to weight system we need to take care
of the relations (\ref{eq:3term}) for {\em external\/} vertices of
\Fd s. This leads us to the following notion.

\begin{dfn}{\rm
The {\em universal enveloping algebra\/} $\cU(L)$ of the $S$-Lie
algebra $L$ is the quotient algebra of the tensor algebra $T^*(L)$ by
the ideal generated by the expressions $a\otimes b - S(a\otimes b) -
f(a,b)$ for $a,b \in L$.
}\end{dfn}
\medskip

Denote by $W_L$ the composition of $F_L$ and the projection $T^*(L)
\rightarrow \cU(L)$.  The following fact makes it possible to pull
$W_L$ down to the space $\cB=\oplus \cB_p$ of \Fd s modulo relations
(\ref{eq:3term}).

\begin{prop} \label{prop:center}
For every $\Gamma \in \cG_{0,*}$ the element $W_L(\Gamma) \in \cU(L)$
is invariant with respect to the canonical action of $L$ on $\cU(L)$.
In other words, $W_L(\Gamma)$ belongs to the $S$-center of $\cU(L)$.
\end{prop}

Now, given a \Fd\ $\Gamma$ and a self-dual Lie $S$-algebra $L$, we can
define $W_L(\Gamma) \in \cU(L)$ to be $W_L(\Gamma')$, where $\Gamma'$
is any element of $\cG_{0,*}$ whose closure (i.e., a \Fd\ obtained by
placing the legs of $\Gamma'$ on the Wilson loop) is
$\Gamma$. The previous proposition guarantees that the result does not
depend on the choice of $\Gamma'$.

Summarizing we come to the following result of \cite{YB}.

\begin{thm}  
For every self-dual Yang--Baxter Lie algebra $L$ there exists a
natural homomorphism of algebras $W_{{L}}: \cA \rightarrow \cU(L)$.
The image of $W_{{L}}$ belongs to the $S$-center $Z_L$ of $\cU(L)$.
\end{thm}
\medskip

The weight system $W_L$ is called {\em the universal weight system}
corresponding to $L$.

\subsection{Examples and particular cases}

\noindent
{\bf 1.~Lie (super)algebras as $S$-Lie algebras.}
\
Perhaps the most familiar example of an $S$-Lie algebra, $L$, is that
of a Lie (super)algebra when $L$ is a ${\Bbb Z}_2$-graded vector space
$L = L_{\bar0} \oplus L_{\bar1}$. Here, the $S$-structure, $S:L\otimes
L\,\to\,L\otimes L,$ is a linear map given by
\begin{displaymath}
S(u\otimes v) = (-1)^{|u| |v|} v \otimes u
\end{displaymath}
for homogeneous elements $u$ and $v$ with degrees $|u|$ and $|v|$,
respectively.  (In the case of an ordinary Lie algebra, $L_{\bar1} =
0$ and $S$ is nothing more than the usual transposition of factors in
the tensor product.)  Relation (3) in the definition of an $S$-Lie
algebra then says that the Lie bracket is $\Z_2$-graded (i.e., the map
$f:L\otimes L \rightarrow L$ preserves the parity). Relations (8) and
(9) become the usual (super-)skew-symmetry and the (super-)Jacobi
identity.  Finally, relations (7), (8) and (5) say, respectively, that
the bilinear form $b$ is invariant, (super-)symmetric and of even
parity with respect to the grading.
\medskip

\noindent
{\bf 2.~Evaluating a diagram (an example).}
\
Consider the following chord diagram $D$:
\medskip
\begin{displaymath}
\Picture{\FullCircle \EndChord[6,0]\EndChord[2,10]\EndChord[4,8]}.
\end{displaymath}
\vspace{20pt}

\noindent
If it is cut open, it becomes the following:
\begin{displaymath}
\begin{picture}(4,1)
\qbezier(0.5,0)(1.5,2)(2.5,0)
\qbezier(2.0,0)(3.0,2)(4.0,0)
\qbezier(3.5,0)(4.5,2)(5.5,0)
\thicklines
\put(0,0){\line(1,0){6}}
\end{picture}
\end{displaymath}
The associated element in $\cU(L)$ is
\begin{displaymath}
W_L(D) = \sum_{a_1,\ldots,a_{10}=1}^{{\rm dim} L} C^{a_1 a_7} C^{a_8
a_9} C^{a_{10} a_6} S^{a_2 a_3}_{a_7 a_8} S^{a_4 a_5}_{a_9 a_{10}} e_{a_1}
e_{a_2} e_{a_3} e_{a_4} e_{a_5} e_{a_6}.
\end{displaymath}
Notice that although we know from Proposition \ref{prop:center} that
$W_L(D)$ belongs to $Z_L$, the precise expression of $W_L(D)$ in terms
of the generators of $Z_L$ requires a computation {\em in $\cU(L)$\/}.
Thus it is a computationally involved matter to determine, for
example, given two diagrams $D$ and $D'$, whether or not $W_L(D)$ and
$W_L(D')$ are the same.  A better method to compute $W_L(D)$ is thus
required, and our recursion formula provides just such a method in the
case where $L$ is $\gl$.
\medskip

\noindent
{\bf 3.~Some special diagrams.}
\  Henceforth, and in order to 
simplify diagrams, when displaying the value of $W_L$
on a chord diagram, we will simply draw the diagram.

\begin{itemize}
\item [(a)] The identity element $1$ in $\cU(L)$ is given by
\begin{displaymath}
\abb
1 = \Picture{\FullCircle}
\end{displaymath}
since $W_{{L}}: \cA \rightarrow \cU(L)$ is a homomorphism of algebras
with unit.

\item [(b)] The Casimir, $c$, (dual to the inner product $b$) is given by
\begin{equation}
\abb
c = \Picture{\FullCircle\EndChord[3,9]}~.
\label{eq:cdiag}
\end{equation}

\item [(c)] The image in $\cU(L)$ of the standard totally $S$-antisymmetric
element in $T^3(L)$ is defined by
\begin{equation}
\abb
y  =
\Picture{
\put(0,0){\line(0,-1){1}}
\put(0,-1){\circle*{0.15}}
\put(0,0){\line(5,3){0.857}}
\put(0,0){\line(-5,3){0.857}}
\put(-0.857,0.514){\circle*{0.15}}
\put(0.857,0.514){\circle*{0.15}}
\FullCircle
}~.\label{eq:ydef}
\end{equation}
\vspace{12pt}

\item [(d)] Finally, the following ``bubble'' diagram $B$ is the image in
$\cU(L)$ of the dual to the Killing form for $L$:
\begin{displaymath}
\abb
B = \Picture{
\Endpoint[3]
\Endpoint[9]
\put(0,0){\circle{0.5}}
\put(0,1){\line(0,-1){0.75}}
\put(0,-1){\line(0,1){0.75}}
\FullCircle
}~,
\end{displaymath}
and can be shown to obey $B = 2y$.
\end{itemize}
\medskip

\noindent
{\bf 4.~The leading order terms.}
\
For $D\in \cD_n$ the leading term of $W_L(D)$ (with respect
to the filtration in $\cU(L)$ inherited from $T(L)$) is equal to 
$c^n$.

In some cases, the subleading term can also be readily
computed.  In those cases for which $B$ is proportional to $c$ (e.g.,
when $L$ a simple Lie algebra) the numerical coefficient of the
$c^{n-1}$ term is proportional to the number of intersections of the
chords in $D$ when the chords are placed in generic position. 
\medskip

\noindent
{\bf 5.~Multiplicativity.}
\
Because $W_{{L}}: \cA \rightarrow \cU(L)$ is a homomorphism of
algebras, we need only consider indecomposable diagrams when 
computing universal weight systems, that is, diagrams whose chords
cannot be split into two non-intersecting subsets.
For example, we see that
\begin{displaymath}
\abb
\Picture{\Chord[1,11]\Chord[3,7]\Chord[5,9]\FullCircle\Endpoint[1]
\Endpoint[3]\Endpoint[5]\Endpoint[7] \Endpoint[9]\Endpoint[11]}
= \Picture{\Endpoint[3]\Endpoint[5]\Endpoint[7]
\Endpoint[9]\Chord[3,7] \Chord[5,9]\FullCircle}\
\Picture{\Endpoint[1]\Endpoint[11]\Chord[1,11]\FullCircle} = c\
\Picture{\Endpoint[3]\Endpoint[7] \Endpoint[9]\Endpoint[5]\Chord[3,7]
\Chord[5,9]\FullCircle}~. 
\end{displaymath}


\section{Representations and invariant tensors of\\ $\gl$}

In this section we present algebraic facts about the Lie superalgebra
$\gl$ necessary for the analysis of the corresponding universal weight
system.  We define $\gl$, examine the subring of its representation
ring generated by cyclic modules, and determine its invariant tensors.
We will also prove an identity between invariant tensors of fourth
order which is crucial in the proof of the recursion relation
(\ref{eq:rec}).

We assume that the ground field $k$ is $\C$.

\subsection{Preliminaries on \gl}

We follow the standard conventions about linear superspaces and Lie
superalgebras which can be found for example in \cite{Berezin,Kac}.

\begin{dfn}{\rm
Let $V \cong \C^{1|1}$ be a $(1|1)$-dimensional vector superspace.
(Recall that an $(m|n)$-dimensional vector superspace, $V$, is a
$\Z_2$-graded vector space $V = V_{\bar 0} \,\oplus V_{\bar 1}$ where
$\dim V_{\bar 0} = m$ and $\dim V_{\bar 1}=n$.)  The Lie superalgebra
of endomorphisms of $V$ is called $\gl$.
}\end{dfn}

The bilinear form 
\begin{equation}
\langle x,y\rangle_{\str} = \str  (xy)
\label{eq:strmetric}
\end{equation}
on $\gl$ is invariant and nondegenerate.  (Recall that the {\em
supertrace\/} of an $(m|n)\times (m|n)$ matrix $M = \pmatrix{A & B\cr
C &  D\cr}$ is defined as   $ \str  M = \tr  A - \tr  D$.)  Therefore,
$\gl$ is a self-dual Lie superalgebra.

Relative to a homogeneous basis $(e_0,e_1)$ for $V$, with $e_0\in
V_{\bar 0}$ and $e_1 \in V_{\bar 1}$, a convenient basis for $\gl$ is
given by the matrices:
\begin{equation}                  \label{eq:basis}
H = \pmatrix{1 & 0\cr 0 & 1\cr},~
G = \pmatrix{0 & 0\cr 0 & 1\cr},~
Q_+ = \pmatrix{0 & 0 \cr 1 & 0 \cr},~
Q_- = \pmatrix{0 & 1 \cr 0 & 0\cr}.
\end{equation}

The even part $L_{\bar 0}$ of $\gl=L=L_{\bar 0}\oplus L_{\bar 1}$ is
spanned by $\{H,G\}$ and the odd part $L_{\bar 1}$ is spanned by
$\{Q_+,Q_-\}$.  The Lie bracket $[-,-]: L^{\otimes 2} \to L$ in this
basis is
\begin{equation}
[G,Q_\pm] = \pm Q_\pm\quad \hbox{and}\quad [Q_+,Q_-]=H,
\label{eq:brackets}
\end{equation}
and zero everywhere else.

The nonvanishing inner products between elements of the basis
(\ref{eq:basis}) are
\begin{equation}
\langle H,G\rangle_{\str} =-1, \quad \langle G,G\rangle_{\str} =-1,
\quad \text{and} \quad 
\langle Q_+,Q_-\rangle_{\str}  = -1~.  \label{eq:strproducts}
\end{equation}

\medskip

The element $H$ belongs to the center of $\gl$ and the quotient Lie
superalgebra $\gl/ \langle H\rangle $ is called $p\gl$.

\subsection{The ring of cyclic modules}

The structure of the representation ring of $\gl$ is not trivial,
since not every finite dimensional module is completely reducible.
The situation for the cyclic modules (modules generated by one
element) is much simpler and luckily this is all we need for our
purposes.

We will only consider those modules $V$ in which $G$ and $H$ act
diagonally.  Since $H$ is in the center of $\gl$, it must be a scalar
operator in any cyclic module: $H v = {\lambda}\,v$ for all $v\in V$.
Let $V = \bigoplus_{\gamma} V_{\gamma}$ be the decomposition of $V$
into eigenspaces of $G$.  We will call any ${\gamma}$ for which
$V_{\gamma}\neq 0$ a {\em weight\/} of $V$.  We will now study the
cyclic modules.

Because $Q_\pm^2 = \frac{1}{2} [Q_\pm,Q_\pm] = 0$, any cyclic module
is at most $(2|2)$-dimensional.  Indeed, let $v$ be a cyclic vector.
Then the module is spanned by the following four elements: $v$, $Q_+
v$, $Q_- v$, and $Q_+ Q_- v$, which need not all be different from
zero.  We can distinguish the following types of cyclic modules, which
are depicted in the Table \ref{tab:cyc}.

\begin{table}[h!]
\begin{displaymath}
\setlength{\unitlength}{1pt}
\begin{array}{c@{\hspace{40pt}}c@{\hspace{40pt}}c}
\I_{\gamma}\quad
\begin{picture}(40,20)(-20,0)
\put(0,0){\circle*{2}}
\put(4,3){\makebox(0,0){${}^{\gamma}$}}
\end{picture}
&
\II^+_{\gamma}\quad
\begin{picture}(40,20)(-20,0)
\put(-12,4){\makebox(0,0){${}^{\gamma-1}$}}
\put(-10,0){\circle*{2}}
\put(10,0){\circle*{2}}
\put(15,4){\makebox(0,0){${}^{\gamma}$}}
\put(-8,0){\vector(1,0){16}}
\end{picture}
&
\II^-_{\gamma}\quad
\begin{picture}(40,20)(-20,0)
\put(-12,4){\makebox(0,0){${}^{\gamma-1}$}}
\put(-10,0){\circle*{2}}
\put(10,0){\circle*{2}}
\put(15,4){\makebox(0,0){${}^{\gamma}$}}
\put(8,0){\vector(-1,0){16}}
\end{picture}
\\ [20pt]
\II_{{\lambda},{\gamma}}\quad
\begin{picture}(40,10)(-20,0)
\put(-13,5){\makebox(0,0){${}^{{\lambda},{\gamma-1}}$}}
\put(-10,0){\circle*{2}}
\put(10,0){\circle*{2}}
\put(15,5){\makebox(0,0){${}^{{\lambda},{\gamma}}$}}
\put(-8,2){\vector(1,0){16}}
\put(8,-2){\vector(-1,0){16}}
\end{picture}
&
\III^0_{\gamma}\quad
\begin{picture}(40,10)(-20,0)
\put(-21,4){\makebox(0,0){${}^{\gamma-2}$}}
\put(-20,0){\circle*{2}}
\put(2,4){\makebox(0,0){${}^{\gamma-1}$}}
\put(0,0){\circle*{2}}
\put(20,0){\circle*{2}}
\put(25,4){\makebox(0,0){${}^{\gamma}$}}
\put(2,0){\vector(1,0){16}}
\put(-2,0){\vector(-1,0){16}}
\end{picture}
&
\III_{\gamma}\quad
\begin{picture}(40,10)(-20,0)
\put(-22,5){\makebox(0,0){${}^{\gamma-2}$}}
\put(-20,0){\circle*{2}}
\put(1,14){\makebox(0,0){${}^{\gamma-1}$}}
\put(0,10){\circle*{2}}
\put(0,-10){\circle*{2}}
\put(20,0){\circle*{2}}
\put(25,4){\makebox(0,0){${}^{\gamma}$}}
\put(-18,2){\vector(2,1){16}}
\put(-2,-8){\vector(-2,1){16}}
\put(18,2){\vector(-2,1){16}}
\put(2,-8){\vector(2,1){16}}
\end{picture}
\\ [20pt]
\end{array}
\end{displaymath}
\caption{Indecomposable cyclic modules.  Right arrows denote 
the action of $Q_+$ and left arrows that of $Q_-$.}
\label{tab:cyc}
\end{table}

\def\np{\vspace{12pt}\noindent}
\np
\underline{\em Type $\I$\/}~~Both $Q_+  v$ and $Q_- v$ vanish.\\
In this case ${\lambda}=0$, since ${\lambda}v=Hv=[Q_+,Q_-]v=0$.  If $G
v = {\gamma}\,v$, then we say that the module is of type
$\I_{\gamma}$.  Such modules are irreducible.  They are either
$(0|1)$- or $(1|0)$-dimensional, depending on the parity of $v$.

\np
\underline{\em Type $\II$\/}~~{\em Only\/} one of $Q_-  v$ and $Q_+ v$
is zero {\em and\/} ${\lambda}\neq 0$.\\ 
If, say $Q_- v \neq 0$, then
the module $V$ is spanned by $v$ and $w = Q_- v$ or, equivalently, by
$w$ and $Q_+ w$.  Both $v$ and $w$ are cyclic vectors, but $v$ in this
case is the highest-weight vector.  If $v$ has weight ${\gamma}$, then
we say that the module is of type $\II_{{\lambda},{\gamma}}$.  These
modules are $(1|1)$-dimensional and also irreducible.

\np
\underline{\em Type $\II^+$\/}~~In this case $Q_-  v = 0$ and $w
= Q_+  v \neq 0$ {\em but\/} ${\lambda}=0$.\\
Because of this condition, $Q_- w = 0$.  If $v$ has weight ${\gamma}$
we say that the module is of type $\II^+_{\gamma}$.  These modules are
$(1|1)$-dimensional, reducible but indecomposable.

\np
\underline{\em Type $\II^-$\/}~~Now $Q_+  v = 0$ and $w = Q_-  v \neq
0$ and ${\lambda}$ is still zero.\\
In this case the highest-weight vector is not the cyclic vector. If
$w$ has weight ${\gamma}$ (so that $v$ has weight ${\gamma}-1$) we say
that the module is of type $\II^-_{\gamma}$.  Again such a module is
$(1|1)$-dimensional, reducible and indecomposable.

\np
\underline{\em Type $\III$\/}~~Both $Q_+  v$ and $Q_-  v$ are 
nonzero.\\
Here we have to distinguish two cases depending on whether ${\lambda}$
is zero or not. If ${\lambda}\neq 0$, then $Q_+ Q_- v$ is also
nonzero.
 
In this case we say that the module $V$ is of type
$\III_{{\lambda},{\gamma}}$, where ${\gamma}$ is the weight of the
highest weight vector $Q_+ v$.  Such modules are $(2|2)$-dimensional
and reducible; but they are also decomposable:
$\III_{{\lambda},{\gamma}} \cong \II_{{\lambda},{\gamma}} \oplus
\II_{{\lambda},{\gamma}-1}$.

When ${\lambda}=0$, the vector $Q_+ Q_- v$ coincides with $ - Q_- Q_+
v$, but this vector might vanish.  If it does not, then $V$ is a
$(2|2)$-dimensional reducible module which, if the highest-weight
vector $Q_+ v$ has weight ${\gamma}$, we denote by $\III_{\gamma}$.
It is indecomposable.

Finally, if $Q_+ Q_- v = 0$, then we have a $(1|2)$- or
$(2|1)$-dimensional reducible indecomposable module which, if it has
highest-weight ${\gamma}$, is denoted by $\III^0_{\gamma}$.

\np
Note that the defining representation is of type $\II_{1,1}$ whereas
the adjoint representation is of type $\III_1$. The Lie superalgebra
$p\gl$ as a $\gl$-module is of type $\III_1^0$.

\medskip

The following is a summary of our case study.

\begin{prop} \label{prop:cyclic}
Each diagonalizable indecomposable cyclic $\gl$-module $V$ is
isomorphic to one of the modules
\begin{equation}
\label{eq:cyclic}
\I_{\gamma},\quad \II_{{\lambda},{\gamma}} \ ({\lambda}\neq 0), \quad
\II^\pm_{\gamma}, \quad \III^0_{\gamma}, \quad \III_{\gamma}.
\end{equation}
\end{prop}
\medskip

Notice that there is a two-fold ambiguity implicit in our notation.
Every module type described above comes in two flavors depending on
the choice of parity of the highest-weight vector.  When it is
necessary to specify this parity, we will do it explicitly, for
example $\III_{\gamma}^{\text{odd}}$, etc.

Direct sums of the cyclic modules considered above form a subring of
the representation ring of $\gl$, as the next proposition shows.

\begin{prop}
The above cyclic modules have the following multiplication
table under tensor product.

\begin{eqnarray}
\I_\gamma \otimes \I_{\gamma'} &\cong&
\I_{\gamma+\gamma'}\nonumber\\
\I_\gamma \otimes \II_{\lambda,\gamma'} &\cong&
                     \II_{\lambda,\gamma+\gamma'}\nonumber\\
\I_\gamma \otimes \II^\pm_{\gamma'} &\cong&
                  \II^\pm_{\gamma+\gamma'}\nonumber\\
\I_\gamma \otimes \III^0_{\gamma'} &\cong&
                \III^0_{\gamma+\gamma'}\nonumber\\
\I_\gamma \otimes \III_{\gamma'} &\cong&
\III_{\gamma+\gamma'}\nonumber\\
\II_{\lambda,\gamma} \otimes \II_{\lambda',\gamma'} &\cong&
\cases{\II_{\lambda+\lambda',\gamma+\gamma'} \oplus
\II_{\lambda+\lambda',\gamma+\gamma'-1}&
              for $\lambda+\lambda'\neq 0$,\cr 
 \III_{\gamma+\gamma'}& for $\lambda =
-\lambda'$.}\nonumber\\
\II_{\lambda,\gamma} \otimes \II^\pm_{\gamma'} &\cong&
\II_{\lambda,\gamma+\gamma'} \oplus
\II_{\lambda,\gamma+\gamma'-1}\nonumber\\
\II^\pm_\gamma \otimes \II^\pm_{\gamma'} &\cong&
\II^\pm_{\gamma+\gamma'} \oplus
 \II^\pm_{\gamma+\gamma'-1}\nonumber\\
\II^+_\gamma \otimes \II^-_{\gamma'} &\cong&
\III_{\gamma+\gamma'}\label{eq:III}\\
\II_{\lambda,\gamma} \otimes \III^0_{\gamma'} &\cong&
\II_{\lambda,\gamma+\gamma'} \oplus
 \II_{\lambda,\gamma+\gamma'-1} \oplus
\II_{\lambda,\gamma+\gamma'-2}\nonumber\\
\II^+_\gamma \otimes \III^0_{\gamma'} &\cong&
\II^+_{\gamma+\gamma'} \oplus
 \III_{\gamma+\gamma'-1}\nonumber\\
\II^-_\gamma \otimes \III^0_{\gamma'} &\cong&
\III_{\gamma+\gamma'} \oplus
 \II^-_{\gamma+\gamma'-2}\nonumber\\
\III^0_\gamma \otimes \III^0_{\gamma'} &\cong&
\II^+_{\gamma+\gamma'} \oplus
 \III^0_{\gamma+\gamma'-1} \oplus 2 \I_{\gamma+\gamma'-2}
\oplus \II^-_{\gamma+\gamma'-3}\nonumber
\\
\II_{\lambda,\gamma} \otimes \III_{\gamma'} &\cong&
\II_{\lambda,\gamma+\gamma'} \oplus 2\,
            \II_{\lambda,\gamma+\gamma'-1}  
 \oplus \II_{\lambda,\gamma+\gamma'-2}\nonumber\\
\II^\pm_\gamma \otimes \III_{\gamma'} &\cong&
\III_{\gamma+\gamma'} \oplus
\III_{\gamma+\gamma'-1}\nonumber\\ 
\III^0_\gamma \otimes \III_{\gamma'} &\cong&
\III_{\gamma+\gamma'} \oplus 2
\III_{\gamma+\gamma'-1}\nonumber\\ 
\III_\gamma \otimes \III_{\gamma'} &\cong&
\III_{\gamma+\gamma'} \oplus 2\, \III_{\gamma+\gamma'-1}
\oplus 
 \III_{\gamma+\gamma'-2}\label{eq:tensor}
\end{eqnarray}
\end{prop}
\begproof
Let $V$ and $V'$ be two cyclic modules with highest vectors $v$ and
$v'$, weights $\gamma$ and $\gamma'$, and $H$ acting by scalars
$\lambda$ and $\lambda'$ respectively. The module $W=V\otimes V'$ has
highest vector $w=v\otimes v'$ with $Hw = (\lambda + \lambda') w$ and
$Gw = (\gamma + \gamma') w$.  Together with the classification of
cyclic modules before Proposition \ref{prop:cyclic}, this gives the
first nine statements of our multiplication table.

The products of $\II_{\lambda,\gamma}$ and modules with $\lambda = 0$
will have $\lambda \neq 0$ and, therefore, will be sums of modules of
the same type.

When $V=\II^-_\gamma$ and $V'=\III^0_{\gamma'}$ with cyclic vectors
$v$ and $v'$ respectively, the vectors $v\otimes v'$ and $v \otimes
Q_- v'$ generate two non-intersecting cyclic submodules in $V\otimes
V'$ of types $\III_{\gamma+\gamma'}$ and $\II^-_{\gamma+\gamma'-2}$
respectively.  Therefore, $\II^-_\gamma \otimes \III^0_{\gamma'} \cong
\III_{\gamma+\gamma'} \oplus \II^-_{\gamma+\gamma'-2}$.

The cases $\II^+_\gamma \otimes \III^0_{\gamma'}$ and $\III^0_\gamma
\otimes \III^0_{\gamma'}$ are handled similarly.

The last four formulas follow from the previous ones using
$\III_\gamma \cong \II^+_{0}\otimes\II^-_\gamma$ (\ref{eq:III}) and
the associativity of the tensor product.
\eproof

\subsection{Invariant tensors}

Let $L$ denote the adjoint module of $\gl$. The tensor powers
$L^{\otimes n}$ have a natural $\gl$-module structure.  Since $L$ is
of type $\III_1$, it follows from (\ref{eq:tensor}) that $L^{\otimes
n}$ decomposes into a direct sum of modules of types $\III_p$ for
integer $p$.  In fact, iterating (\ref{eq:tensor}) we obtain

\begin{prop}
 \begin{equation}
  L^{\otimes n} \cong \III_1^{\otimes n} \cong \bigoplus_{\ell =
  -n+2}^{n}\, {2n-2 \choose n-\ell}\, \III_\ell \label{eq:decomp}
 \end{equation}
\end{prop}
\eproof
\medskip

Now we will compute the number of linearly independent invariant
tensors of order $n$, i.e., the dimension of the subspace

\begin{displaymath}
\Inv (L^{\otimes n}) = \{ v \in L^{\otimes n}\ |\ gv=0 \quad
\text{for\ all\ } g \in \gl\}~.
\end{displaymath}

Nonzero invariants in cyclic modules can exist only when
$\lambda{=}0$, and the following lemma tells us for which values of
$\gamma$ they appear in the indecomposable cyclic modules
(\ref{eq:cyclic}).

\begin{lemma}\label{invariants}
Let $V$ be one of the cyclic modules (\ref{eq:cyclic}).  Then the
dimension of the space $\Inv V $ is either $0$ or $1$ and $\dim \Inv V
= 1$ if and only if $\lambda=0$ and
\begin{eqnarray*}
\gamma &=& \cases{0&for $V$ of type $\I_\gamma$ or $\II^{+}_\gamma$\cr
1&for $V$ of type $\II^{-}_\gamma$ or $\III_\gamma$ \cr
0 \hbox{ or } 2 & for $V$ of type $\III^{0}_\gamma$}.
\end{eqnarray*}
\end{lemma}
\eproof

It follows that $\dim\Inv (L^{\otimes n})$ is equal to the
multiplicity of $\III_1$ in the direct sum decomposition of
$L^{\otimes n}$, which from (\ref{eq:decomp}) equals ${2n-2 \choose
n-1}$.  In summary,

\begin{prop}\label{tensinv}
The number of linearly independent $\gl$-invariant tensors of order
$n$ on $L$ is given by ${2n-2 \choose n-1}$.
\end{prop}
\eproof

For $n=1$ there is only one independent invariant tensor $H$.  For
$n=2$ there are two\footnote{In order to avoid notational clutter, we
will hereafter omit the $\otimes$ from the notation for $\gl$-tensors.
Therefore we will understand $x\,y$ to mean $x\otimes y$ and $x^2 =
x\otimes x$, etc.}: one is $H^2$, and the other one is given by
\begin{equation}
C = G\,H + H\,G + Q_- \,Q_+ - Q_+ \,Q_-~.\label{eq:dualmetric}
\end{equation}
Notice that they are both (super)symmetric.  Together they span the
space of $\gl$-invariants in $S^2 L$.  In particular, $C$ is
nondegenerate and therefore induces a $\gl$-isomorphism $L \cong L^*$.
This isomorphism allows us to identify $S^2 L$ with the space of
invariant symmetric bilinear forms on $L$.  Under this identification,
we associate with the invariant tensor $\alpha C + \beta H^2$ the
bilinear form
\begin{equation}
\langle Q_+,Q_-\rangle = \alpha~,~\langle H,G\rangle = 
\alpha~\hbox{and}~\langle G,G\rangle = \beta. \label{eq:metric}
\end{equation}
Provided $\alpha\neq 0$ this bilinear form is nondegenerate; that is,
it defines an invariant metric on $\gl$.  Therefore there are several
ways to make $\gl$ a self-dual Lie superalgebra.  The metric
(\ref{eq:metric}) coincides with $\langle -, - \rangle_{\str}$ when
$\alpha = -1$ and $\beta = -1$.

Weight systems $W_L$ corresponding to different choices of $\alpha\neq
0~\hbox{and}~\beta$ in (\ref{eq:metric}) are equivalent, since the
transformation
\begin{displaymath}
H\mapsto \alpha^{-1} H,~G\mapsto G -\frac{\beta}{2\alpha}H,~Q_\pm 
\mapsto \frac{1}{\sqrt{\alpha}}Q_\pm
\end{displaymath}
is an automorphism of $\gl$ that takes the bilinear form
(\ref{eq:metric}) to the metric dual to $C$ (it corresponds to the
choice $\alpha=1~\hbox{and}~\beta=0$). From now on we will use this
metric on $\gl$ in our computations of the universal weight system
$W_{gl(1|1)}$.

\subsection{Tensor diagrammar}

The language of Feynman graphs of subsection \ref{sec:FG} provides
a convenient tool to treat (invariant) tensors in Lie superalgebras
graphically.  We will use such pictures for denoting components of
tensors as well.  Let $L$ be any vector superspace with a fixed
homogeneous basis $\{e_i\}$.  Relative to this basis, every linear map
$\varphi: L^{\otimes m}
\to L^{\otimes n}$ is determined by the $(\dim L)^{m+n}$ numbers
$(\varphi_{i_1 i_2\cdots i_m}^{j_1 j_2 \cdots j_n})$
defined by
\begin{displaymath}
\varphi(e_{i_1}\, e_{i_2} \cdots e_{i_m}) = 
\varphi_{i_1 i_2\cdots i_m}^{j_1 j_2 \cdots j_n}\,
e_{j_1}\,e_{j_2}\cdots e_{j_n}~,
\end{displaymath}
where we use the summation convention (and all the $\otimes$ between
$e_i$ are omitted).  We will represent it
graphically by the following diagram:
\begin{displaymath}
\begin{picture}(5,3)
\put(1.5,0){\line(0,1){1}}
\put(1.25,0.5){\makebox(0,0)[l]{$\scriptstyle 1$}}
\put(2,0){\line(0,1){1}}
\put(1.75,0.5){\makebox(0,0)[l]{$\scriptstyle 2$}}
\put(2.75,0.5){\makebox(0,0){$\cdots$}}
\put(3.5,0){\line(0,1){1}}
\put(3.6,0.5){\makebox(0,0)[l]{$\scriptstyle n$}}
\put(1,1){\framebox(3,1){$\varphi$}}
\put(1.5,2){\line(0,1){1}}
\put(1.25,2.5){\makebox(0,0)[l]{$\scriptstyle 1$}}
\put(2,2){\line(0,1){1}}
\put(1.75,2.5){\makebox(0,0)[l]{$\scriptstyle 2$}}
\put(2.75,2.5){\makebox(0,0){$\cdots$}}
\put(3.5,2){\line(0,1){1}}
\put(3.6,2.5){\makebox(0,0)[l]{$\scriptstyle m$}}
\thicklines
\put(0,0){\line(1,0){5}}
\put(0,3){\line(1,0){5}}
\end{picture}~.
\nonumber
\end{displaymath}

For  an arbitrary vector superspace $L$, there are two canonical maps:
the identity and the symmetry.  The identity map ${\mathrm id}: L\to L$,
with components $\delta_i^j$, is represented by
\begin{displaymath}
{\delta_i}^j \leftrightarrow
\begin{picture}(2,2)
\put(1,0){\line(0,1){2}}
\put(1,-0.3){\makebox(0,0){${}_j$}}
\put(1,2.3){\makebox(0,0){${}^i$}}
\thicklines
\put(0,0){\line(1,0){2}}
\put(0,2){\line(1,0){2}}
\end{picture}~;
\nonumber
\end{displaymath}
\vspace{12pt}

\noindent
whereas the symmetry $S: L^{\otimes 2} \to L^{\otimes 2}$ has
components $S_{ij}^{k\ell} = (-1)^{|i||j|} \delta_i^\ell
\delta_j^k$, where $|i|$ is the parity of $e_i$; and is represented
by

\vspace{10pt}
\begin{displaymath}
S_{ij}^{k\ell} \leftrightarrow
\begin{picture}(3,2)
\put(1,0){\line(1,2){1}}
\put(2,0){\line(-1,2){1}}
\put(1,-0.3){\makebox(0,0){${}_k$}}
\put(2,-0.3){\makebox(0,0){${}_\ell$}}
\put(1,2.3){\makebox(0,0){${}^i$}}
\put(2,2.3){\makebox(0,0){${}^j$}}
\thicklines
\put(0,0){\line(1,0){3}}
\put(0,2){\line(1,0){3}}
\end{picture}~.
\nonumber
\end{displaymath}

Now let $L$ be a Lie superalgebra.  The Lie bracket is a linear map
$L^{\otimes 2} \to L$, or in component form $[e_i,e_j] = f_{ij}^k
e_k$. It will be depicted as
\begin{equation}
f_{ij}^k \leftrightarrow
\begin{picture}(4,2.5)
\qbezier(1,2)(2,0)(3,2)
\put(2,0){\line(0,1){1}}
\put(2,-0.3){\makebox(0,0){${}_k$}}
\put(1,2.3){\makebox(0,0){${}^i$}}
\put(3,2.3){\makebox(0,0){${}^j$}}
\thicklines
\put(0,0){\line(1,0){4}}
\put(0,2){\line(1,0){4}}
\end{picture}~.
\label{eq:f-ijk}
\end{equation}

If, in addition, $L$ is self-dual, there is a canonical quadratic
invariant tensor.  If we let $C_{ij} = \langle e_i,e_j \rangle$ denote
the coefficients of the invariant metric in this basis,\footnote{In
the previous section, the metric was denoted by $b$.}{} the following
is an invariant tensor in $L^{\otimes 2}$:

\begin{equation}
C \equiv C^{ij} e_i \, e_j,
\nonumber
\end{equation}
where $(C^{ij})$ is the matrix inverse of $(C_{ij})$; that is, $C^{ij}
C_{jk} = \delta^i_k$.  The diagram which represents this invariant
tensor is
\begin{equation}
C^{ij} \leftrightarrow
\begin{picture}(4,1)
\qbezier(1,0)(2,2)(3,0)
\thicklines
\put(0,0){\line(1,0){4}}
\put(1,-0.3){\makebox(0,0){${}_i$}}
\put(3,-0.3){\makebox(0,0){${}_j$}}
\end{picture}~.
\label{eq:quadinv}
\end{equation}
The invariant metric itself can be represented diagrammatically as
follows:
\smallskip

\begin{displaymath}
C_{ij} \leftrightarrow 
\begin{picture}(4,1)
\qbezier(1,1)(2,-1)(3,1)
\thicklines
\put(0,1){\line(1,0){4}}
\put(1,1.3){\makebox(0,0){${}^i$}}
\put(3,1.3){\makebox(0,0){${}^j$}}
\end{picture}~.
\end{displaymath}

Just as in the case of Feynman graphs considered in Section 1, there
are two ways to combine diagrams of invariant tensors together to make
a third diagram: tensor product and composition.  Because the
identity, the symmetry, the Lie bracket, the metric and the quadratic
tensor are $L$-invariant, so will be any diagram obtained by gluing
these together.  For example, gluing two copies of (\ref{eq:quadinv})
and (\ref{eq:f-ijk}) we can obtain a cubic invariant.  Concretely, we
define the invariant tensor $F\in L^{\otimes 3}$ by the diagram
\def\CubicInv{
\begin{picture}(4,1)
\qbezier(1,0)(2,2)(3,0)
\put(2,0){\line(0,1){1}}
\thicklines
\put(0,0){\line(1,0){4}}
\end{picture}
}
\begin{eqnarray}
\CubicInv &=& \begin{picture}(4,3)
\qbezier(0.5,1)(1,0)(1.5,1)
\put(1,0){\line(0,1){0.5}}
\put(2,0){\line(0,1){1}}
\put(3,0){\line(0,1){2}}
\qbezier(1.5,1)(1.75,1.5)(2,1)
\put(0,2){\line(1,0){4}}
\put(0.5,1){\line(0,1){1}}
\qbezier(0.5,2)(1.75,4)(3,2)
\put(0,1){\line(1,0){4}}
\thicklines
\put(0,0){\line(1,0){4}}
\end{picture}\nonumber\\
\ \nonumber \\
&=& f_{\ell m}^i C^{m j} C^{\ell k} e_i \, e_j \, e_k~.
\label{eq:cubinv}
\end{eqnarray}

Although many invariant tensors can be constructed from $C^{ij}$ and
$f_{ij}^k$ by gluing, it is important to keep in mind that not every
invariant tensor is of this form.  Indeed, returning now to $\gl$ we
see this immediately.

Already for order one, the only linearly-independent $\gl$-invariant
tensor, $H$, cannot be written in terms of the quadratic
(\ref{eq:quadinv}) and cubic (\ref{eq:cubinv}) invariants.  If we wish
to be able to depict $H$ graphically we must introduce a new symbol:
\begin{displaymath}
H \leftrightarrow
\begin{picture}(2,2)
\put(1,0){\line(0,1){1.5}}
\put(1,1.5){\circle*{0.3}}
\thicklines
\put(0,0){\line(1,0){2}}
\end{picture}
\nonumber
\end{displaymath}

For $n{=}2$ we saw that there are two linearly independent invariant
tensors: $H^2$ and the quadratic invariant $C$ (\ref{eq:quadinv}),
which in the case of $\gl$ is given by (\ref{eq:dualmetric}).  But
these are not all the invariant tensors of order two that we can
draw. Take, for example, the ``bubble'' $B$, defined by the diagram
\def\Bubble{
\begin{picture}(4,2)
\put(1.5,0){\oval(2,2)[tl]}
\put(2.5,0){\oval(2,2)[tr]}
\put(2,1){\circle{1}}
\thicklines
\put(0,0){\line(1,0){4}}
\end{picture}}
\begin{displaymath}
B = \Bubble~.
\nonumber
\end{displaymath}
This is clearly an invariant tensor, since it is constructed out of
gluing invariant tensors (in fact, it is dual to the Killing form)
\begin{eqnarray*}
\Bubble &=&
\begin{picture}(4,3)
\put(1,0){\line(0,1){.5}}
\put(3,0){\line(0,1){.5}}
\qbezier(0.5,1)(1,0)(1.5,1)
\qbezier(2.5,1)(3,0)(3.5,1)
\put(0,1){\line(1,0){4}}
\qbezier(1.5,1)(2,2)(2.5,1)
\put(0.5,1){\line(0,1){1}}
\put(3.5,1){\line(0,1){1}}
\put(0,2){\line(1,0){4}}
\qbezier(0.5,2)(2,4)(3.5,2)
\thicklines
\put(0,0){\line(1,0){4}}
\end{picture}
\\
\ \\
&=& C^{kn} C^{\ell m} f_{k\ell}^i f_{mn}^j\, e_i \, e_j~;
\end{eqnarray*}
but since according to Proposition \ref{tensinv} there are only two
linearly independent quadratic invariant tensors, there must be a
linear relation between $H^2$, $C$ and $B$. Indeed, bubbles
burst:

\begin{prop}\label{prop:bubble}
\be
\Bubble =
-2~ \begin{picture}(4,2)
\put(1.5,0){\oval(2,2)[tl]}
\put(2.5,0){\oval(2,2)[tr]}
\put(1.5,1){\circle*{0.3}}
\put(2.5,1){\circle*{0.3}}
\thicklines
\put(0,0){\line(1,0){4}}
\end{picture}
\qquad \text{or \ } B = - 2 H^2.
\ee
\end{prop}
\begproof
From the fact that the Killing form of $p\gl$ vanishes, or
equivalently from the fact that $G$ is not a linear combination of
commutators (\ref{eq:brackets}), it follows that $B$ is proportional
to $H^2$.  It is then a simple matter of finding the constant of
proportionality.
\eproof

According to Proposition \ref{tensinv} there are six linearly
independent invariant tensors of order three.  We can write five of
them using the invariant tensors $H$, $C$ and $F$ that we have
introduced before:
\begin{eqnarray*}
\begin{picture}(4,1)
\put(1,0){\line(0,1){1}}
\put(2,0){\line(0,1){1}}
\put(3,0){\line(0,1){1}}
\put(1,1){\circle*{0.3}}
\put(2,1){\circle*{0.3}}
\put(3,1){\circle*{0.3}}
\thicklines
\put(0,0){\line(1,0){4}}
\end{picture} &=& H^3\\
\begin{picture}(4,1)
\put(3,0){\line(0,1){1}}
\put(3,1){\circle*{0.3}}
\qbezier(0.5,0)(1.5,2)(2.5,0)
\thicklines
\put(0,0){\line(1,0){4}}
\end{picture} &=& C \, H\\
\begin{picture}(4,1)
\put(2,0){\line(0,1){1}}
\put(2,1){\circle*{0.3}}
\qbezier(1,0)(5,2)(3,0)
\thicklines
\put(0,0){\line(1,0){4}}
\end{picture} &=& C^{ij} e_i\,H\,e_j\\
\begin{picture}(4,1)
\put(1,0){\line(0,1){1}}
\put(1,1){\circle*{0.3}}
\qbezier(1.5,0)(2.5,2)(3.5,0)
\thicklines
\put(0,0){\line(1,0){4}}
\end{picture} &=& H \, C\\
\CubicInv &=& -H \, Q_- \, Q_+ - H \, Q_+ \, Q_- - Q_- \, H \,
Q_+ + Q_- \, Q_+ \, H\\ 
&&{}- Q_+ \, H \, Q_- + Q_+ \, Q_- \, H~;\\
\end{eqnarray*}
but the sixth invariant is of a different kind:
\begin{eqnarray*}
T &=& G^2 \, H + G \, H \, G  + H \, G^2 + G \, Q_- \, Q_+ - G \, Q_+
\, Q_- +  Q_- \, G \, Q_+\\
&&{}+ Q_- \, Q_+ \, G - Q_+ \, G \, Q_- - Q_+ \,
Q_- \, G .
\end{eqnarray*}

We will not attempt to write down twenty linearly independent
invariant tensors that exist in $L^{\otimes 4}$.  We simply note that
five tensors in $\Inv(L^{\otimes 4})$ which can be constructed from
the quadratic and cubic tensors satisfy a linear relation.

\subsection{A fundamental relation}

The following result is crucial in the proof of the recursion formula
(\ref{eq:rec}).

\begin{thm}
Let $K$ be the invariant tensor corresponding to the following
diagram:
\def\DiagramK{
\begin{picture}(5,2)
\qbezier(0.5,0)(1.25,1.5)(2,0)
\qbezier(3,0)(3.75,1.5)(4.5,0)
\qbezier(1.25,0.75)(2.5,2)(3.75,0.75)
\thicklines
\put(0,0){\line(1,0){5}}
\end{picture}}
\begin{eqnarray*}
\DiagramK &=&
\begin{picture}(6,2)
\put(0.5,0){\line(0,1){1}}
\put(2,0){\line(0,1){0.5}}
\put(4,0){\line(0,1){0.5}}
\put(5.5,0){\line(0,1){1}}
\put(0,1){\line(1,0){6}}
\qbezier(0.5,1)(1,2)(1.5,1)
\qbezier(1.5,1)(2,0)(2.5,1)
\qbezier(2.5,1)(3,2)(3.5,1)
\qbezier(3.5,1)(4,0)(4.5,1)
\qbezier(4.5,1)(5,2)(5.5,1)
\thicklines
\put(0,0){\line(1,0){6}}
\end{picture}\\
\ \\
&=& C^{im} C^{np} C^{q\ell} f_{mn}^j f_{pq}^k e_i \, e_j
\, e_k \, e_\ell
\end{eqnarray*}
and let
\begin{eqnarray}
M &=&
\begin{picture}(3,1)
\put(0.75,0){\oval(1,1)[tl]}
\put(1.25,0){\oval(1,1)[tr]}
\put(1,0.5){\circle{0.5}}
\put(2,0){\oval(1.5,1)[t]}
\thicklines
\put(0,0){\line(1,0){3}}
\end{picture} + 
\begin{picture}(3,1)
\put(1,0){\oval(1.5,1)[t]}
\put(1.75,0){\oval(1,1)[tl]}
\put(2.25,0){\oval(1,1)[tr]}
\put(2,0.5){\circle{0.5}}
\thicklines
\put(0,0){\line(1,0){3}}
\end{picture} - 
\begin{picture}(3,1.5)
\put(1.25,0){\oval(1,1)[tl]}
\put(1.75,0){\oval(1,1)[tr]}
\put(1.5,0.5){\circle{0.5}}
\put(1.5,0){\oval(2.5,2)[t]}
\thicklines
\put(0,0){\line(1,0){3}}
\end{picture} -
\begin{picture}(3,1.5)
\put(1.5,0){\oval(1.5,1)[t]}
\put(1.25,0){\oval(2,2)[tl]}
\put(1.75,0){\oval(2,2)[tr]}
\put(1.5,1){\circle{0.5}}
\thicklines
\put(0,0){\line(1,0){3}}
\end{picture}\nonumber\\
&=& B_{13}C_{24} + C_{13} B_{24} - C_{14} B_{23} -
B_{14}C_{23}~,\nonumber
\end{eqnarray}
then 
\begin{equation}
K = \frac{1}{2} M
\label{eq:tensoreq}
\end{equation}
(where $C_{24}B_{13} = C^{j\ell}B^{ik}\, e_i \, e_j\,
e_k\, e_\ell$, etc.)
\end{thm}

\begproof
Notice that both $K$ and $M$ belong to $\Inv(S^2 (\bigwedge^2 L))
\subset L^{\otimes 4}$. 

Let $L_0$ be the quotient superalgebra $p\gl = \gl/ \langle H \rangle$
considered as $\gl$-module, and

\begin{displaymath}
p : \Inv(S^2 (\bigwedge\nolimits^2 L)) \ra \Inv(S^2
(\bigwedge\nolimits^2 L_0))
\end{displaymath}
the map induced by the projection $L \ra L_0$.  From Lemma \ref{s2l2}
below, it follows that the space $\Inv(S^2 (\bigwedge^2 L))$ is
two-dimensional and, in fact, spanned by the two tensors
\begin{eqnarray*}
M &=& B_{13}C_{24} + C_{13} B_{24} - C_{14} B_{23} - B_{14}C_{23}\\
\text{and\  \ } N &=&
\begin{picture}(4,1)
\qbezier(0.5,0)(1.5,2)(2.5,0)
\qbezier(1.5,0)(2.5,2)(3.5,0)
\thicklines
\put(0,0){\line(1,0){4}}
\end{picture} - 
\begin{picture}(4,1.5)
\qbezier(0.5,0)(2,2.5)(3.5,0)
\qbezier(1.5,0)(2,1)(2.5,0)
\thicklines
\put(0,0){\line(1,0){4}}
\end{picture}\\
&=& C_{13}C_{24} - C_{14} C_{23}~, 
\end{eqnarray*}
such that $M$ is a generator of $\Ker(p)$ and $p(N)$ spans
$\Inv(S^2 (\bigwedge^2 L_0))$.

On the other hand, the invariant tensor $K$ is given by
\begin{eqnarray*}
K &=& - H \, Q_- \, H \, Q_+ + H \, Q_- \, Q_+ \, H + H \, Q_+ \, H \,
Q_- - H \, Q_+ \, Q_- \, H\\
&&{}+ Q_- \, H^2 \, Q_+ - Q_- \, H \, Q_+\, H - Q_+ \, H^2 \, Q_- +
Q_+ \, H \, Q_-\, H~.
\end{eqnarray*}
Therefore $K$ belongs to $\Ker(p)$ and has to be proportional to $M$.
By comparing the coefficients at, say, $Q_- H^2 Q_+$ in $K$ and $M$ we
find that the constant of proportionality is equal to ${1\over 2}$.
\eproof

It remains to prove the lemma used above.

\begin{lemma}\label{s2l2}
\begin{eqnarray*}
(i) & &\qquad \dim \Inv (S^2(\bigwedge\nolimits^2 L))  =  2;\\
(ii) & &\qquad \dim \Inv (S^2(\bigwedge\nolimits^2 L_0)) = 1;\\
(iii) & &\qquad \dim \Ker\left(\Inv (S^2(\bigwedge\nolimits^2 L))
\ra \Inv (S^2(\bigwedge\nolimits^2 L_0)\right)  =  1.
\end{eqnarray*}
\end{lemma}
\begproof \

Paying closer attention to the calculation behind
(\ref{eq:tensor}), it follows that
\begin{eqnarray*}
&&\bigwedge\nolimits^2 \III_g^{\mathrm odd} \cong \III_{2g}^{\mathrm 
even} \oplus \III_{2g-2}^{\mathrm even}\quad\hbox{and}\quad
S^2 \III_g^{\mathrm odd} \cong 2 \III_{2g-1}^{\mathrm odd}\\
&&\bigwedge\nolimits^2 \III_g^{\mathrm even} \cong 2
\III_{2g-1}^{\mathrm odd}\quad\hbox{and}\quad 
S^2 \III_g^{\mathrm even} \cong \III_{2g}^{\mathrm even} \oplus
\III_{2g-2}^{\mathrm even}~.
\end{eqnarray*}
Since the adjoint module is of type $\III_1^{\mathrm odd}$, it 
follows
that
\begin{eqnarray*}
S^2\bigwedge\nolimits^2 \III_1^{\mathrm odd} &\cong&
S^2\left(\III_{2}^{\mathrm even} \oplus \III_{0}^{\mathrm
even}\right)\\
&\cong& \III_4^{\mathrm even} \oplus 2 \III_2^{\mathrm even} \oplus 2
\III_1^{\mathrm odd} \oplus 2 \III_0^{\mathrm even} \oplus
\III_{-2}^{\mathrm even}~.
\end{eqnarray*}
Part (i) of the lemma now follows from Lemma \ref{invariants}.

The case (ii) is proved similarly by using the fact that the module
$L_0$ is of type $\III_1^0$.

The statement (iii) follows from (i), (ii) and the fact that $p$ is
surjective (since $p(N) \neq 0$).  By Proposition \ref{prop:bubble}
the ``bubble'' vanishes in $L_0^{\otimes 2}$, therefore $M \in \Ker
(p)$.
\eproof
\vspace{1cm}

\begin{rem}
{\rm
A similar (and even slightly simpler) argument can be used to
prove the key relation between invariant $sl_2$-tensors needed for
computation of the universal $sl_2$-weight system in \cite{CV}.  In
this case $\dim \Inv (S^2(\bigwedge\nolimits^2 L)) = 1$.  }
\end{rem}
\bigskip


\section{$U(gl(1|1))$-valued weight systems}

\subsection{Fundamental relation on the level of weight systems}

Now we will rewrite the relation (\ref{eq:tensoreq}) in terms of the
universal $\cU(\gl)$-valued \ws\ $W$.  The weight system $W$ is
defined as in Example 1.3.2.  Take a Feynman diagram $F$ and cut open
the Wilson loop anywhere but at an external leg.  This gives us an
invariant tensor, which is not unique, since it depends on where we
cut the diagram open.  However, as it was explained in Section
\ref{sec:FG}, the image $W(F)$ of this tensor in the universal
enveloping algebra $\cU(\gl)$ is well-defined and by
Proposition \ref{prop:center} belongs to the center of $\cU(\gl)$.

It is known (cf. \cite{Berezin}) that the center of $\cU(\gl)$ is
the polynomial algebra $\C\,[h,c]$ where $h$ and $c$ are the images in
$\cU(\gl)$ of the invariant tensors $H$ and $C$, respectively.
However since Feynman diagrams only possess trivalent vertices and
propagators, the image of $W$ is not all of $\C\,[h,c]$ but the
subalgebra $\C\,[c,y]$ where $y$ was defined in equation
(\ref{eq:ydef}).  This fact is an immediate corollary from the
recursion relation (\ref{eq:rec}).  With our choice of invariant
metric, it is easy to check that
\begin{equation}
y = -h^2~.\label{eq:ydeftoo}
\end{equation}

It follows from the previous discussion that any linear relation in
the tensor algebra like the one in (\ref{eq:tensoreq}) can be inserted
in any Feynman diagram to yield a linear relation for the
corresponding weight system.  Indeed, in any fixed chord diagram of
order $n{-}2$, we can insert (\ref{eq:tensoreq}) to obtain the
following relation for the universal weight system $W$:

\begin{equation}
\Picture{
\DottedCircle
\FullChord[5,8]
\FullChord[2,11]
\put(-0.55,-0.183){\line(3,1){1.1}}
} = - {1\over 2}\left(
\Picture{
\DottedCircle
\FullChord[2,5]
\put(0.05,-0.35){\circle{0.4}}
\qbezier(-0.5,-0.866)(-0.4,-0.5)(-0.163,-0.427)
\qbezier(0.866,-0.5)(0.65,-0.2)(0.273,-0.3)
\Endpoint[11]
\Endpoint[8]
\Arc[8]
\Arc[11]
}
+ 
\Picture{
\DottedCircle
\FullChord[8,11]
\put(-0.05,0.35){\circle{0.4}}
\qbezier(0.5,0.866)(0.4,0.5)(0.163,0.427)
\qbezier(-0.866,0.5)(-0.65,0.2)(-0.273,0.3)
\Endpoint[2]
\Endpoint[5]
\Arc[2]
\Arc[5]
} - \Picture{
\DottedCircle
\Endpoint[2]
\Endpoint[8]
\qbezier(0.5,0.866)(0.4,0.7)(0.3,0.533)
\put(0.197,0.362){\circle{0.4}}
\put(-0.5,-0.866){\line(3,5){0.625}}
\FullChord[5,11]
\Arc[2]
\Arc[8]
}
- 
\Picture{
\DottedCircle
\Endpoint[5]
\Endpoint[11]
\qbezier(0.866,-0.5)(0.7,-0.4)(0.533,-0.3)
\put(0.362,-0.197){\circle{0.4}}
\put(-0.866,0.5){\line(5,-3){1.042}}
\FullChord[2,8]
\Arc[5]
\Arc[11]
}
\right)~,\label{eq:crucial}
\end{equation}
where the $n{-}2$ original chords are not pictured but are the same
for all five diagrams.

The following identity for $W$ is an immediate consequence of
Proposition \ref{prop:bubble} and equation (\ref{eq:ydeftoo}):
\begin{equation}
\Picture{
\DottedCircle
\Endpoint[3]
\Endpoint[9]
\put(0,0){\circle{0.5}}
\put(0,1){\line(0,-1){0.75}}
\put(0,-1){\line(0,1){0.75}}
\Arc[3]
\Arc[9]
}
= 
2y\,\Picture{
\DottedCircle
\Arc[3]
\Arc[9]
}~.\label{eq:bubblesoff}
\end{equation}
\vspace{10pt}

Finally, as a direct corollary of the relations (\ref{eq:crucial}) and
(\ref{eq:bubblesoff}) we have:
\begin{prop}
The universal $\cU(\gl)$-valued weight system $W$ satisfies the
following relation:
\begin{displaymath}
\Picture{\DottedCircle\FullChord[5,8]\FullChord[2,11]\put(-0.55,-0.183){\line(3,1){1.1}}}
= y\ \left(
\Picture{\DottedCircle\Arc[2]\Arc[8]\FullChord[5,11]} +
\Picture{\DottedCircle\Arc[5]\Arc[11]\FullChord[2,8]} -
\Picture{\DottedCircle\Arc[8]\Arc[11]\FullChord[2,5]} -
\Picture{\DottedCircle\Arc[2]\Arc[5]\FullChord[8,11]}\right)~.
\end{displaymath}
\label{cor:fund}
\end{prop}
\eproof

\subsection{The recursion formula}

We now prove the recursion formula stated in the Introduction, but
we shall first need a lemma.

\begin{lemma} (The Eight-Term Relations) \label{lem:8term}
Let $R_1,\ldots,R_4$ be defined by the following formulas:
\begin{eqnarray*}
R_1 &=& \Picture{\DottedCircle\FullChord[1,4]\FullChord[5,7]
\FullChord[8,10]} - 
\Picture{\DottedCircle\FullChord[1,4]\FullChord[5,8]
\FullChord[7,10]} -
\Picture{\DottedCircle\FullChord[1,5]\FullChord[4,7]
\FullChord[8,10]} +
\Picture{\DottedCircle\FullChord[1,5]\FullChord[4,8]
\FullChord[7,10]} \\ [20pt]
R_2 &=& \Picture{\DottedCircle\FullChord[1,8]\FullChord[4,11]
\FullChord[5,7]} -  
\Picture{\DottedCircle\FullChord[1,7]\FullChord[4,11]
\FullChord[5,8]} -
\Picture{\DottedCircle\FullChord[1,8]\FullChord[5,11]
\FullChord[4,7]} +
\Picture{\DottedCircle\FullChord[1,7]\FullChord[4,8]
\FullChord[5,11]} \\ [20pt]
R_3 &=& \Picture{\DottedCircle\FullChord[1,7]\FullChord[2,4]
\FullChord[8,10]} -
\Picture{\DottedCircle\FullChord[1,8]\FullChord[2,4]
\FullChord[7,10]} -
\Picture{\DottedCircle\FullChord[1,4]\FullChord[2,7]
\FullChord[8,10]} +
\Picture{\DottedCircle\FullChord[1,4]\FullChord[2,8]
\FullChord[7,10]} \\ [20pt]
R_4 &=& \Picture{\DottedCircle\FullChord[1,4]\FullChord[5,11]
\FullChord[7,10]} -
\Picture{\DottedCircle\FullChord[1,4]\FullChord[5,10]
\FullChord[7,11]} -
\Picture{\DottedCircle\FullChord[1,5]\FullChord[4,11]
\FullChord[7,10]} +
\Picture{\DottedCircle\FullChord[1,5]\FullChord[4,10]
\FullChord[7,11]}~,\\[10pt]
\end{eqnarray*}
and let $S_1,\ldots,S_4$ be defined by:
\begin{eqnarray*}
S_1 &=&
\Picture{\DottedCircle\FullChord[1,7]\Arc[4]\Arc[5]\Arc[8]\Arc[10]} +
\Picture{\DottedCircle\FullChord[4,10]\Arc[1]\Arc[5]\Arc[7]\Arc[8]} -
\Picture{\DottedCircle\FullChord[4,7]\Arc[1]\Arc[5]\Arc[8]\Arc[10]} -
\Picture{\DottedCircle\FullChord[1,10]\Arc[4]\Arc[5]\Arc[7]\Arc[8]}\\
[20pt]
S_2 &=&
\Picture{\DottedCircle\FullChord[8,10]\Arc[1]\Arc[4]\Arc[5]\Arc[7]} +
\Picture{\DottedCircle\FullChord[2,4]\Arc[5]\Arc[7]\Arc[8]\Arc[11]} -
\Picture{\DottedCircle\FullChord[4,7]\Arc[1]\Arc[5]\Arc[8]\Arc[11]} -
\Picture{\DottedCircle\FullChord[1,11]\Arc[4]\Arc[5]\Arc[7]\Arc[8]}\\
\end{eqnarray*}

\clearpage

\begin{eqnarray*}
S_3 &=&
\Picture{\DottedCircle\FullChord[1,8]\Arc[2]\Arc[4]\Arc[7]\Arc[10]} +
\Picture{\DottedCircle\FullChord[4,10]\Arc[1]\Arc[2]\Arc[7]\Arc[8]} -
\Picture{\DottedCircle\FullChord[4,7]\Arc[1]\Arc[2]\Arc[8]\Arc[10]} -
\Picture{\DottedCircle\FullChord[1,10]\Arc[2]\Arc[4]\Arc[7]\Arc[8]}\\
[20pt]
S_4 &=&
\Picture{\DottedCircle\FullChord[1,7]\Arc[4]\Arc[5]\Arc[10]\Arc[11]} +
\Picture{\DottedCircle\FullChord[4,10]\Arc[1]\Arc[5]\Arc[7]\Arc[11]} -
\Picture{\DottedCircle\FullChord[4,7]\Arc[1]\Arc[5]\Arc[10]\Arc[11]} -
\Picture{\DottedCircle\FullChord[1,10]\Arc[4]\Arc[5]\Arc[7]\Arc[11]}~.
\\[10pt]
\end{eqnarray*}
Then for all $i=1,\ldots,4$,
\begin{displaymath}
R_i = y S_i~.
\end{displaymath}
\end{lemma}

\begproof
The left hand sides of the first, third, and fourth equations are obtained
by applying the three-term relations to the diagram
\begin{displaymath}
\Picture{
\DottedCircle
\FullChord[2,4]
\FullChord[8,10]
\put(0,-0.7){\line(0,1){1.4}}
}~.
\end{displaymath}
\vspace{10pt}

On the other hand, applying Proposition \ref{cor:fund} to the same
diagram gives rise to the right hand sides.  The left hand side of the
second equation is obtained by applying the transposition operator,
$S$, to the right two ends of the preceding diagram:
\begin{displaymath}
\Picture{
\DottedCircle
\FullChord[4,11]
\FullChord[1,8]
\put(-0.2,-0.55){\line(0,1){1.1}}
}~.
\end{displaymath}
\vspace{10pt}

Doing the same to the right hand sides of Proposition \ref{cor:fund}
gives rise to its right hand side.
\eproof

\begin{thm}\label{thm:rec}
Let $D$ be a chord diagram, ``$a$'' a fixed chord in $D$, and $b_1$,
$b_2$, \ldots be all the chords of $D$ intersecting $a$. Denote by
$D_a$ (resp.  $D_{a,i}$, $D_{a,ij}$) the diagram $D - a$ (resp. $D - a
- b_i$, $D- a - b_i - b_j$). Then
\begin{eqnarray*}
W(D) &=& c W(D_a) - y\, \sum_i W(D_{a,i}) \\
&&{} +  y\, \sum_{i<j} \left(W\left(D_{a,ij}^{+-}\right) +
W\left(D_{a,ij}^{-+}\right)- W\left(D_{a,ij}^l\right) -
W\left(D_{a,ij}^r\right)\right)~,
\end{eqnarray*}
where $D_{a,ij}^{+-}$ (resp. $D_{a,ij}^{-+}$; $D_{a,ij}^l$;
$D_{a,ij}^r$) is the diagram obtained by adding to $D_{a,ij}$ a new
chord connecting the left end of $b_i$ and the right end of $b_j$
(resp. the right end of $b_i$ and the left end of $b_j$; the left ends
of $b_i$ and $b_j$; the right ends of $b_i$ and $b_j$) assuming that
the chord $a$ is drawn vertically. See equation (\ref{eq:pixD}) for a
pictorial description of the diagrams appearing in the recursion
relation.
\end{thm}

\begproof
The proof is a double induction similar to the one done by
Chmutov-Varchenko \cite{CV} for $sl_2$.  Let $D$ be a chord diagram of
order $\ord D=n$ and let ``$a$'' be a chord in $D$.  Draw $D$ so that
$a$ is a vertical line where the number $\ell = \ell(D,a)$ of chord
endpoints to the left of $a$ is less than or equal to the number of
chord endpoints to the right of $a$.  Assume that $\ell\geq 2$.  We
take as the induction hypothesis that the recursion formula holds for
all pairs $(E,b)$ consisting of a chord diagram $E$ and a chord $b$ in
$E$, such that $\ord E\leq n$ and $\lambda(E,b)<\ell$.  For the pair
$(D,a)$ described above, there are seven possible configurations of
chords which have endpoints immediately to the left of $a$:
\vspace{20pt}
\begin{displaymath}
\def\VerticalChord{
\put(0.5,0.866){\line(0,-1){1.732}}
\Endpoint[2]\Arc[2]
\Endpoint[10]\Arc[10]
 }
\begin{array}{ccccccc}
\Picture{\DottedCircle\FullChord[3,9]
\FullChord[4,1] \FullChord[8,11]}
&\Picture{\DottedCircle\VerticalChord
\FullChord[3,5] \FullChord[7,9]}
&\Picture{\DottedCircle\FullChord[3,9]
\FullChord[4,11]\FullChord[8,1]}
&\Picture{\DottedCircle\VerticalChord
\FullChord[3,7] \FullChord[9,5]}
&\Picture{\DottedCircle\FullChord[3,9]
\FullChord[4,6] \FullChord[8,0]}
&\Picture{\DottedCircle\FullChord[3,9]
\FullChord[4,0]\FullChord[8,6]}
&\Picture{\DottedCircle\FullChord[3,9]
\FullChord[4,8]}
\end{array}
\end{displaymath}
\vspace{20pt}

The vertical line is the chord $a$ in each case. In the first six
cases, we denote by $u$ and $v$ respectively the upper and lower
chords adjacent to $a$ in the diagram. We will work through the
induction step for a diagram, $D$, that looks like the first
configuration. Configurations two to six follow similarly. (The
induction step for the last configuration is trivial because a chord
intersects $a$ if and only if it intersects the adjacent chord.) Apply
the first eight-term relation of Lemma \ref{lem:8term} to $D$ to
obtain:  
\begin{eqnarray*}
\Picture{\DottedCircle\FullChord[3,9] \FullChord[4,1]\FullChord[8,11]}
={} &&y\ \left(
\Picture{\DottedCircle\FullChord[1,8]\Arc[3]\Arc[4]\Arc[9]\Arc[11]} +
\Picture{\DottedCircle\FullChord[4,11]\Arc[1]\Arc[3]\Arc[8]\Arc[9]} -
\Picture{\DottedCircle\FullChord[4,8]\Arc[1]\Arc[3]\Arc[9]\Arc[11]} -
\Picture{\DottedCircle\FullChord[1,11]\Arc[3]\Arc[4]\Arc[8]\Arc[9]}
\right)\\[20pt]
&& {}+
\Picture{\DottedCircle\FullChord[1,3]\FullChord[4,9]\FullChord[8,11]}
+
\Picture{\DottedCircle\FullChord[1,4]\FullChord[3,8]\FullChord[9,11]}
-
\Picture{\DottedCircle\FullChord[1,3]\FullChord[4,8]\FullChord[9,11]}
\\[10pt]
\end{eqnarray*}

Let us denote the terms of the right hand side of this equation by
$X_1$, $X_2$, \ldots, $X_7$, respectively. Since the terms on the
right hand side have fewer endpoints to the left of $a$ than $D$, by
the induction hypothesis, we obtain the following:
\begin{enumerate}
\item $X_1 = y\, W\left(D_{a,uv}^{-+}\right)$
\item $X_2 = y\, W\left(D_{a,uv}^{+-}\right)$
\item $X_3 = -y\, W\left(D_{a,uv}^{l}\right)$
\item $X_4 = -y\, W\left(D_{a,uv}^{r}\right)$
\item $X_5 = c  W\left(D_a\right) - y\,\left( \sum_i
W\left(D_{a,i}\right) +  W\left(D_{a,v}\right) - \sum_{i<j}
W(\Lambda_{a,ij})\right.\\ \phantom{X_5 = {}} -  \sum_i
W\left(\Lambda_{a,iv}\right)\biggr)$
\item $X_6 =  c  W\left(D_a\right) -y\,\left( \sum_i
W\left(D_{a,i}\right) + W\left(D_{a,u}\right) - \sum_{i<j}
W(\Lambda_{a,ij})\right.\\ \phantom{X_6 = {}} - \left. \sum_j
W\left(\Lambda_{a,uj}\right)\right)$
\item $X_7 =  -c  W\left(D_a\right) + y\,\left( \sum_i
W\left(D_{a,i}\right)  - \sum_{i<j} W(\Lambda_{a,ij})\right)$
\end{enumerate}
where $\Lambda_{a,ij}$ is given by
\begin{equation}
\Lambda_{a,ij} = D_{a,ij}^{+-} + D_{a,ij}^{-+} - D_{a,ij}^{l} -
D_{a,ij}^{r}~.\label{eq:Lambda}
\end{equation}
Adding all the terms $X_1, \ldots, X_7,$ we obtain the right hand side
of the recursion formula. This completes the induction step.

It remains to show the induction base: $\ell = 0$ and $\ell=1$.
If $\ell=0$, then no chord in $D$ intersects $a$, so that $D$ is
a connected sum of the chord diagram $\Theta$ of order $1$
and a chord diagram $D_a$ of order $n-1$. Therefore, $W(D) =
c\,W(D_a)$ which agrees with the recursion formula.  If $\ell=1$,
there is only one one chord $i$, say, intersecting $a$:
\begin{displaymath}
\Picture{
\DottedCircle\FullChord[5,7]\FullChord[0,6]
\put(-0.55,0.3){\makebox(0,0){${}^a$}}
\put(-0.2,0.1){\makebox(0,0){${}^i$}}
}
\end{displaymath}
\vspace{10pt}

Using the three-term relation (\ref{eq:3term}),
\begin{eqnarray*}
\Picture{\DottedCircle\FullChord[5,7]\FullChord[0,6]} &=&
\Picture{\DottedCircle\FullChord[5,6]\FullChord[0,7]} -
\Picture{\DottedCircle \put(1,0){\line(-1,0){1.68} }
         \FullChord[5,7]\Arc[0]\Arc[6]\Endpoint[0]   
}\\[20pt]
&=& c\ \Picture{\DottedCircle\Arc[5]\Arc[6]\FullChord[0,7]} -
\half\ 
\Picture{\DottedCircle
     \put(0,0){\circle{0.5}}
     \put(1,0){\line(-1,0){0.75}}     \put(-1,0){\line(1,0){0.75}}
     \Arc[0]\Arc[5]\Arc[7]\Arc[6] \Endpoint[6]\Endpoint[0]
}
= c\ \Picture{\DottedCircle\Arc[5]\Arc[6]\FullChord[0,7]} 
- y\ \Picture{\DottedCircle\Arc[0]\Arc[5]\Arc[6]\Arc[7]}~,
\\[10pt]
\end{eqnarray*}
which is just $c\,W(D_a) - y\,W(D_{a,i})$, in agreement with the
recursion formula.
\eproof

\begin{cor}
The universal weight system $W_{gl(1|1)}$ corresponding to $\gl$ takes
values in the polynomial algebra $\Z\,[c,y]$. Furthermore, the
polynomial corresponding to a chord diagram, $D$, of order $n$ is a
weighted homogeneous polynomial of order $n$ where $c$ and $y$ are
assigned weights $1$ and $2,$ respectively.  The leading term in $c$
of $W_{gl(1|1)}(D)$ is equal to $c^n$, where $n$ is the number of
chords in $D$. The coefficient of the subleading term $c^{n-1} y$ is
the negative of the number of intersections of the chords in $D$ when
the chords are placed in generic position.
\end{cor}

\begproof
The corollary follows from a simple induction on the number of chords.
\eproof

\subsection{Examples of computations}
\label{sec:recursion}

We will now illustrate the use of the recursion formula on a few
examples.  Normalizing $W$ by setting it equal 1 on the diagram with
no chords, the recursion formula yields the value of $W$ on any other
diagram.  Equation (\ref{eq:cdiag}) is an immediate consequence of the
recursion relation.  For the unique indecomposable chord diagram of order 2
we obtain

\begin{eqnarray*}
\Picture{\FullCircle\EndChord[0,6]\EndChord[3,9]} &=& 
c\ \Picture{\FullCircle\EndChord[0,6]} - 
y\ \Picture{\FullCircle}\\[18pt]
&=& c^2 - \, y~.
\end{eqnarray*}

There are two indecomposable chord diagrams of order 3:

\begin{eqnarray*}
\Picture{\FullCircle\EndChord[1,5]\EndChord[7,11] \EndChord[3,9]}
 &=& 
c\ \Picture{\FullCircle\EndChord[1,5]\EndChord[7,11]}  - 
y\ \left(
\Picture{\FullCircle\EndChord[7,11]} + 
\Picture{\FullCircle\EndChord[1,5]}\right)\\ [10pt]
&&{}+y\ \left(
\Picture{\FullCircle\EndChord[5,11]} +
\Picture{\FullCircle\EndChord[1,7]} -
\Picture{\FullCircle\EndChord[5,7]} -
\Picture{\FullCircle\EndChord[1,11]} \right)  \\ [10pt]
&=& c^3 - 2 y c~,
\end{eqnarray*}
and
\begin{eqnarray*}
\Picture{\FullCircle\EndChord[2,8]\EndChord[4,10]\EndChord[0,6]}
 &=& 
c\ \Picture{\FullCircle\EndChord[2,8]\EndChord[4,10]}
- y\ \left(
\Picture{\FullCircle\EndChord[2,8]} + 
\Picture{\FullCircle\EndChord[4,10]}\right)\\ [10pt]
&&{}+ y\ \left(
\Picture{\FullCircle\EndChord[2,4]} +
\Picture{\FullCircle\EndChord[8,10]} -
\Picture{\FullCircle\EndChord[2,10]} -
\Picture{\FullCircle\EndChord[4,8]} \right)  \\ [10pt]
&=& c^3 - 3 y c~.
\end{eqnarray*}

These results are summarized in the following table.

\begin{table}[h!]
\begin{center}
\begin{tabular}{||c|l||c|l||} \hline
$D$ & \multicolumn{1}{c||}{$W_{gl(1|1)}(D)$} & $D$ &
\multicolumn{1}{c||}{$W_{gl(1|1)}(D)$}\\
\hline\hline
&&&\\ [3pt]
\Picture{\FullCircle} & 1 &
\Picture{\FullCircle\EndChord[3,9]}& $c$ \\ [20pt]
\Picture{\FullCircle\EndChord[3,9]\EndChord[0,6]} & $c^2 - y$ &
\Picture{\FullCircle\EndChord[2,6]\EndChord[4,10]\EndChord[0,8]}
&$c^3 - 2 cy$\\ [20pt]
\Picture{\FullCircle\EndChord[2,8]\EndChord[4,10]\EndChord[0,6]}
&$c^3 - 3 cy$ &&\\ [20pt]\hline
\end{tabular}
\end{center}
\caption{Indecomposable chord diagrams of order $\leq 3$}
\end{table}

We will now use these results to perform a more complex computation:

\begin{eqnarray*}
\Picture{\FullCircle\EndChord[1,5]\EndChord[7,11]\EndChord[2,10]\EndChord[4,8]}
&=&
c\ \Picture{\FullCircle\EndChord[1,5]\EndChord[7,11]\EndChord[2,10]} - 
y\, \left(\Picture{\FullCircle\EndChord[7,11]\EndChord[2,10]} +
\Picture{\FullCircle\EndChord[1,5]\EndChord[2,10]}\right)\\
[20pt]
&&{}+ y\ \left(
\Picture{\FullCircle\EndChord[5,11]\EndChord[2,10]} +
\Picture{\FullCircle\EndChord[1,7]\EndChord[2,10]}  -
\Picture{\FullCircle\EndChord[1,11]\EndChord[2,10]} -
\Picture{\FullCircle\EndChord[2,10]\EndChord[4,8]}
\right) \\ [10pt]
&=& c(c^3 - 2cy) - 2 y c^2 = c^4 - 4 c^2 y
\end{eqnarray*}

In the appendix we list the values of the universal weight system
$W_{gl(1|1)}$ for all indecomposable chord diagrams of orders 4 and 5.

\subsection{Deframing the universal weight system}

In this subsection we will prove that the strong weight system
obtained from the universal weight system $W_{gl(1|1)}$ under the
deframing projection $\overline{\cal{W}} \longrightarrow {\cal{W}}$ (see
Section 1.1) is the same as evaluating it at $c=0$.\footnote{The
results of this subsection owe much to conversations with Bill
Spence.}

We start by defining a linear map on the algebra $\cA$ of chord
diagrams $s:\cA \to \cA$ as follows.  If $D$ is a nonempty chord
diagram, then
\begin{displaymath}
s(D) = \sum_{a} D_a
\end{displaymath}
where the sum runs over all chords $a$ in $D$ and where $D_a$ denotes
$D - a$; and if $D$ is the empty chord diagram $s(D)=0$.  We extend
$s$ linearly.  It is clear that $s$ respects the four-term relation,
so it is well defined on $\cA$.

Recall (see Theorem \ref{thm:cdalg}) that $\cA =
\bar{\cA}[\Theta]$ where $\Theta$ is the unique chord diagram of
with one chord.  Let $\psi$ be the restriction to $\bar{\cal{A}}$ of
the map $\lambda: \cal{A} \rightarrow \cal{A}$ defined by
\begin{displaymath}
\lambda(D) = \sum_{k\geq 0} (-1)^k \frac{1}{k!} \Theta^k \cdot
s^k(D)~.
\end{displaymath}
Notice that the sum has only a finite number of nonzero terms when
acting on any chord diagram $D$.

\begin{prop}{\cite{BN,KSA}.}
The deframing projector $\cW \to \overline{\cW}$ is dual to
$\psi$.  In other words, if $w\in \cW_n$ is any weight system, and
$D \in \cA_n$ is a chord diagram, then the deframed weight system
$\overline{w}\in\overline{\cW}_n$ is given by:
\begin{displaymath}
\overline{w}(D) = \psi^*w(D) = w(\psi(D))~.
\end{displaymath}
\end{prop}

\begin{thm}
Let $W=W_{gl(1|1)}$ denote the universal weight system of $gl(1|1)$,
i.e., $W: \cA \to Z(\cU(gl(1|1))) = \C[h,c]$, and let
$\overline{W}=\psi^*(W)$ be its deframing. Then if $D$ is any chord
diagram,
\begin{displaymath}
\overline{W} \in \C[h] \ \mbox{ and }\ \overline{W}(h) = W(h,0)~.
\end{displaymath}
\end{thm}

To prove this we start with a lemma.

\begin{lemma}\label{lem:sisddc}
For any chord diagram $D$,
\begin{displaymath}
W(s(D)) = \frac{\partial}{\partial c} W(D)~.
\end{displaymath}
\end{lemma}

\begproof
We proceed by induction on the order of the chord diagram.  The
induction base is clear for diagrams of order 1.  Assume that the
lemma holds for all chord diagrams of order $<n$.  Let $D$ be a chord
diagram of order $n$, and let $a$ be a chord in $D$.  Then $s(D) = D_a +
\sum_{i\neq a} D_i$.  We will use the recursion relation (\ref{eq:rec})
to compute $W(s(D)) = W(D_a) + \sum_{i\neq a} W(D_i)$.  Since $a$ is a
chord in every $D_i$, for $i\neq a$, we can use the recursion relation
on $a$ to obtain
\begin{displaymath}
W(D_i) = c W(D_{ia}) - y \sum_{j\neq i} W(D_{ia,j}) + y \sum_{j<k
\atop j,k\neq i} W(\Lambda_{ia,jk})~,
\end{displaymath}
where $j$ and $k$ are chords in $D$, and where now
\begin{displaymath}
\Lambda_{ia,jk} = D_{ia,jk}^{+-} + D_{ia,jk}^{-+} - D_{ia,jk}^{l} -
D_{ia,jk}^{r}.
\end{displaymath}
On the other hand, we can also use the recursion relation to compute:
\begin{eqnarray*}
\frac{\partial}{\partial c} W(D) &=& \frac{\partial}{\partial c} 
\left( c W(D_a) - y \sum_j W(D_a,j) + y \sum_{j<k}
W(\Lambda_{a,jk})\right)\\
&=& W(D_a) + c \frac{\partial}{\partial c} W(D_a) - y \sum_j
\frac{\partial}{\partial c} W(D_a,j) +  y \sum_{j<k}
\frac{\partial}{\partial c} W(\Lambda_{a,jk}).
\end{eqnarray*}
By the induction hypothesis, $\frac{\partial}{\partial c} W(D_a) =
W(s(D_a))$ and similarly for $D_{a,j}$ and $\Lambda_{a,jk}$, since all
these diagrams have order $<n$.  Noticing that
\begin{displaymath}
s(D_a) = \sum_i D_{ia}~,\qquad s(D_{a,j}) = \sum_{i\neq j} D_{ia,j}~,
\quad\hbox{and}\quad s(\Lambda_{a,jk}) = \sum_{i\neq j,k}
\Lambda_{ia,jk}~,
\end{displaymath}
and inserting into $\frac{\partial}{\partial c} W(D)$, we obtain that
\begin{displaymath}
\frac{\partial}{\partial c} W(D) = W(s(D))~.
\end{displaymath}
\eproof

\noindent
{\em Proof of the theorem\/}.  Using the lemma we can now compute the
deframed universal weight system $\overline{W}$:
\begin{eqnarray*}
\overline{W}(D) & = & \sum_{k\geq 0} (-1)^k \frac{1}{k!}
W(\Theta^k \cdot s^k D) \\
& = &  \sum_{k\geq 0} (-1)^k \frac{1}{k!} W(\Theta^k) W(s^k D)\\
& = &  \sum_{k\geq 0} (-1)^k \frac{1}{k!} c^k
\frac{\partial^k}{\partial c^k} W(D)~.
\end{eqnarray*}
But since $W(D)$ is a polynomial in $c$, the above expression is
simply the Taylor expansion evaluated at $c=0$.
\eproof

\begin{cor}
With the same notation as is Theorem \ref{thm:rec}, the deframed
universal $\cU(\gl)$-valued weight system $\overline{W}$ takes values
in $\Z[y]$ and obeys the following recursion relation:
\begin{displaymath}
\overline{W}(D) = - y\, \sum_i \overline{W}(D_{a,i}) + y \sum_{i<j}
\overline{W}(\Lambda_{a,ij})~,
\end{displaymath}
where $\Lambda_{a,ij}$ is given by (\ref{eq:Lambda}).
\end{cor}

We will see in the next section that specializing to $y{=}1$, the
deframed universal weight system $\overline{W}$ is precisely the
Alexander--Conway weight system.

\subsection{A bosonic version of $gl(1|1)$}

There is a ``bosonic'' analog of $\gl$.  It is the self-dual
four-dimensional Lie algebra $L$ with a basis $\{\,Q_+, Q_-, H, G\,\}$
consisting of even elements with commutators given by
(\ref{eq:brackets}).  Furthermore, the nonvanishing inner products
between elements of this basis are defined by equations
(\ref{eq:metric}) where $\alpha = 1$ and $\beta = 0$.

\begin{thm}
The universal weight system associated to the Lie algebra $L$ is
isomorphic to the universal weight system associated to the Lie
superalgebra $\gl$.
\end{thm}

\begproof
The proof proceeds similarly  to the case of $\gl$. The only
difference is that in this case  $B = 2 h^2$
rather than the $-2 h^2$ as in the $\gl$ case.  
However, when written in terms of $c$ and $y$, the
fundamental relation (\ref{cor:fund}) still holds.  This can be proven
by direct computation of  the invariant tensors appearing in
this identity.
\eproof
\vspace{.5in}


\section{Relation with the Alexander-Conway polynomial}

Here we will show that deframing the universal $\cU(\gl)$-valued
gives, as it was stated in \cite{YB}, the Alexander-Conway weight system.
This yields further evidence of the supersymmetric nature of the
Alexander knot invariant (cf. \cite{KS,RS1}).

\subsection{The Alexander-Conway polynomial}

\begin{dfn} \label{dfn:conw}
{\rm
The {\em Alexander-Conway\/} polynomial 
\begin{displaymath}
\nabla(K)= c_0(K)+c_1(K)z+\ldots c_p(K)z^p
\end{displaymath}
is a $\Z[z]$-valued link invariant uniquely defined by the following
properties:
\begin{enumerate}
\item $\nabla(\text{unknot})=1$
\item $\nabla(K_+)-\nabla(K_-) = z\,\nabla(K_{\mid\mid})$, where
$K_+$, $K_-$, and $K_{\mid\mid}$ are three oriented links that differ
only inside a small ball as indicated:
\begin{displaymath}
\def\ACKII{
\begin{picture}(2,2)(-1,-1)
\put(0,0){\circle{2}}
\qbezier(-0.707,-0.707)(0,0)(-0.707,0.707)
\qbezier(0.707,-0.707)(0,0)(0.707,0.707)
\put(-0.607,0.607){\vector(-1,1){0.1414}}
\put(0.607,0.607){\vector(1,1){0.1414}}
\end{picture}}
\def\ACKPlus{
\begin{picture}(2,2)(-1,-1)
\put(0,0){\circle{2}}
\put(-0.707,-0.707){\vector(1,1){1.414}}
\put(0.707,-0.707){\line(-1,1){0.6}}
\put(-0.107,0.107){\vector(-1,1){0.6}}
\end{picture}}
\def\ACKMinus{
\begin{picture}(2,2)(-1,-1)
\put(0,0){\circle{2}}
\put(-0.707,-0.707){\line(1,1){0.6}}
\put(0.107,0.107){\vector(1,1){0.6}}
\put(0.707,-0.707){\vector(-1,1){1.414}}
\end{picture}}
\mathop{\ACKPlus}_{K_+} \qquad \mathop{\ACKMinus}_{K_-} \qquad
\mathop{\ACKII}_{K_{\mid\mid}}~.
\end{displaymath}
\end{enumerate}
}
\end{dfn}
\medskip

From \ref{dfn:conw}.2 it follows that
\begin{equation}
\nabla(L_1\amalg L_2)=0, 
\label{eq:2comp}
\end{equation}
where $L_1$ and $L_2$ are two links in $\R^3$ separated by a plane.
\medskip

If $K_0, K_+$, and $K_-$ are singular knots (or links) as in
(\ref{eq:vasrel}), then
\begin{equation}
\nabla(K_0)=\nabla(K_+)-\nabla(K_-)=z\,\nabla(K_{\mid\mid})~.
\label{eq:convas}
\end{equation}

By induction we see that for a singular knot $K$ with $n$
self-intersections $\nabla(K)$ is divisible by $z^n$.  Therefore, the
coefficient $c_n$ of ${\nabla}(K)$ at $z^n$ is a Vassiliev invariant
of order $\le n$.

The corresponding weight systems are easy to compute.  For a chord
diagram $D$ of order $n$ denote by $\nabla(D)$ the value of $c_n$ on
$D$.  As before, we will drop the $\nabla(\cdots)$ from the diagrams
in order to simplify the notation.

First notice a simple graphical interpretation of the
Conway--Vassiliev skein relation (\ref{eq:convas}).
\begin{equation}
\Picture{\DottedCircle\Arc[3]\Arc[9]\Chord[3,9]}
=
\Picture{
{\thicklines
\qbezier[25](0,1)(-0.12,1)(-0.259,0.966)
}
\qbezier[4](-0.259,0.966)(-0.518,0.897)(-0.707,0.707)
\qbezier[4](-0.707,0.707)(-0.897,0.518)(-0.966,0.259)
\qbezier[4](-0.966,0.259)(-1.04,0)(-0.966,-0.259)
\qbezier[4](-0.966,-0.259)(-0.897,-0.518)(-0.707,-0.707)
\qbezier[4](-0.707,-0.707)(-0.518,-0.897)(-0.259,-0.966)
{\thicklines
\qbezier[25](-0.259,-0.966)(-0.12,-1)(0,-1)
\put(0,-1){\line(0,1){2}}
}}  
  \hspace{-35pt}
\Picture{
{\thicklines
\qbezier[25](0,1)(0.12,1)(0.259,0.966)
}
\qbezier[4](0.707,-0.707)(0.897,-0.518)(0.966,-0.259)
\qbezier[4](0.966,-0.259)(1.04,0)(0.966,0.259)
\qbezier[4](0.966,0.259)(0.897,0.518)(0.707,0.707)
\qbezier[4](0.707,0.707)(0.518,0.897)(0.259,0.966)
\qbezier[4](0.259,-0.966)(0.518,-0.897)(0.707,-0.707)
{\thicklines
\qbezier[25](0,-1)(0.12,-1)(0.259,-0.966)
\put(0,-1){\line(0,1){2}}
}
}
\label{eq:band}
\\[20pt]
\end{equation}

\begin{rem}
{\rm Here we are dealing with links, therefore we may have diagrams
with more than one Wilson loop; and as before, we assume that diagrams
may contain other chords connecting the dotted arcs, provided that
they are the same in all diagrams appearing in the same relation.  }
\end{rem}

Now, for a chord diagram $D$ we construct an oriented compact surface
with boundary $\S_D$ as follows.  Take a sphere $\S_{\bf O}$ with a
hole (or holes) whose boundary is identified with the Wilson loop
(loops) of $D$.  Replace the chords by narrow non-intersecting ribbons
and attach them to the boundary of $\S_{\bf O}$. For example, for
\begin{displaymath}
\abb
D = \Picture{\FullCircle\Chord[0,6]\Chord[2,10]\Chord[4,8]} \quad
\hbox{we get $S_D =$ \parbox{1in}{\epsfbox{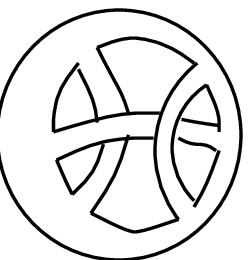}}}.
\end{displaymath}

From formulas (\ref{eq:2comp}) and (\ref{eq:band}) we obtain the
following rule for computing $\nabla(D)$.
\begin{prop}      \label{prop:AC}
The value of the \ws\ $\nabla$ on a chord diagram $D$ is equal to $1$
if $\S_D$ has only one boundary component and $0$ otherwise.
\end{prop}
\ \eproof

For example, for the diagram $D$ above, $\S_D$ has $2$ boundary
components, therefore $\nabla(D)=0$.
\medskip

\subsection{Relation with the universal weight system}

As we saw in Section 3.4, deframing of the universal $gl(1|1)$ weight
system $W_{gl(1|1)}$ is the same as specializing it at $c=0$. The
following theorem shows that it coincides with the Alexander-Conway
weight system.

\begin{thm} \label{th:conw}
The Alexander-Conway weight system $\nabla$ coincides with the
specialization of the universal $gl(1|1)$ weight system $ W_{gl(1|1)}$
for $c=0$ and $y=1$.
\end{thm}

\begproof
Let $W_0$ denote the specialization of $W_{gl(1|1)}$ for $c=0$ and
$y=1$.  Since the condition $W_{gl(1|1)}({\large \bf O})=1$ and the
fundamental relation (\ref{cor:fund}) completely determine
$W_{gl(1|1)}$, it will be enough to verify that the Alexander--Conway
weight system $\nabla$ satisfies the relation (\ref{cor:fund}) for
$c=0, y=1$.

Thus, the theorem follows from the following proposition. 
\eproof

\begin{prop}    \label{prop:rec}
The Alexander--Conway weight system $\nabla$ satisfies the relation
\begin{equation}
\Picture{\DottedCircle\Arc[5]\Arc[7]\Chord[5,7]\Arc[1]\Arc[11]\Chord[1,11]
\put(-0.65,0){\line(1,0){1.33}}}
= \Picture{\DottedCircle\Arc[1]\Arc[7]\Arc[5]\Arc[11]\Chord[5,11]} +
\Picture{\DottedCircle\Arc[5]\Arc[11]\Arc[1]\Arc[7]\Chord[1,7]} -
\Picture{\DottedCircle\Arc[7]\Arc[11]\Arc[1]\Arc[5]\Chord[1,5]} -
\Picture{\DottedCircle\Arc[1]\Arc[5]\Arc[7]\Arc[11]\Chord[7,11]}~.
\label{eq:conwrec}
\end{equation}
\end{prop}
\ \medskip

The proposition follows from several lemmas. 

\begin{lemma}\label{lem:2}
Let $D$ be a chord diagram and $d$ one of its chords.  
Then $\nabla(D) = 1$ if the surface $\S_{D-d}$ has precisely two
boundary components and the endpoints of $d$ belong to different
components and ${\nabla}(D)=0$ otherwise.
\end{lemma}
\begproof  
If the endpoints of $d$ belong to the same boundary component $c$ of
$\S_{D-d}$, then by adding a band corresponding to $d$ we cut $c$ into
two pieces thus increasing the number of components by one.  If $d$
connects two different components, then the band fuses them into one.

The statement now follows from Proposition \ref{prop:AC}.
\eproof

\begin{lemma}\label{lem:3}
Let
\begin{displaymath}
\abb
F = \Picture{
\DottedCircle
\put(0,0){\line(0,-1){1}}
\put(0,0){\line(5,3){0.857}}
\put(0,0){\line(-5,3){0.857}}
\Arc[1]\Arc[5]\Arc[9]
}
\end{displaymath}
\vspace{20pt}

\noindent
be a Feynman diagram with only one trivalent vertex whose edges form
the $Y$ graph (so that $D=F-Y$ is a chord diagram).
Then  $\nabla(F)= \mp \nabla(D)$, namely
\begin{eqnarray*}
{\nabla}(F) &=& {\cases{ 
0 &  \ if $\S_D$ has more than one boundary component, \cr
-1 & \ if $\partial(\S_D)$ has one component and the cyclic order on
the \cr
 & \  legs of $Y$ agrees with the orientation of $c$, \cr
1 & \ otherwise.}}
 \end{eqnarray*}
\end{lemma}

\begproof
Denote by $a$, $b$ and $c$ the boundary components of $\Sigma_D$ to
which the legs of $Y$ are attached. 
We have 
\begin{displaymath}
\nabla(F)= \nabla(D_1) - \nabla(D_2),
\end{displaymath}
where
\begin{displaymath}
D_1 =
\Picture{\DottedCircle
\qbezier(0.866,0.5)(0.1,0)(0.087,-0.996)
\qbezier(-0.866,0.5)(-0.1,0)(-0.087,-0.996)
\Arc[1]\Arc[5]\Arc[9]
}
\qquad \text{and}\qquad
D_2 =
\Picture{\DottedCircle
\qbezier(0.866,0.5)(-0.11,0)(-0.087,-0.996)
\qbezier(-0.866,0.5)(0.11,0)(0.087,-0.996)
\Arc[1]\Arc[5]\Arc[9]
}~.
\end{displaymath}
\\[12pt]

If $\partial\Sigma_D$ has other components besides $a,b$ and $c$, then
$\nabla(F)=0$ by Lemma \ref{lem:2}, since both $D_1$ and $D_2$ have
more than one boundary component.

If, say, $a\ne b$ and $a \ne c$, then $\nabla(D_1)=\nabla(D_2)$ and
$\nabla(F)=0$ since the chord connecting $a$ and $b$ fuses these
components into one without affecting $c$.

At last, if $a=b=c$ and the cyclic order of the endpoints of $Y$
agrees with the orientation of this component, then
\begin{displaymath}
\nabla(F)= \nabla(D_1) - \nabla(D_2)= 0 - 1 =-1,
\end{displaymath}
 and $\nabla(F)=1$ otherwise.
\eproof
\medskip
 
\begin{lemma} \label{lem:4}
Let 
$$ F= 
 \Picture{\DottedCircle\Arc[5]\Arc[7]\Chord[5,7]\Arc[1]\Arc[11]\Chord[1,11]
 \put(-0.66,0){\line(1,0){1.34}}
 \put(-1.25,0.68){\text{\scriptsize d}} 
 \put(-1.25,-0.77){\text{\scriptsize a}} 
 \put(0.70,0.68){\text{\scriptsize c}} 
 \put(0.75,-0.77){\text{\scriptsize b}} 
    } \                     
$$
\\[20pt]

\noindent
be a \Fd\ with exactly
two trivalent vertices connected by a 
propagator, i.e.,
$F=D + K$, where $K$ is the graph 
\Picture{\Chord[5,7]\Chord[1,11]
 \put(-0.66,0){\line(1,0){1.34}}
 \put(-0.99,0.6){\text{\scriptsize D}} 
 \put(-0.99,-0.72){\text{\scriptsize A}} 
 \put(0.73,0.6){\text{\scriptsize C}} 
 \put(0.73,-0.75){\text{\scriptsize B}} 
    } \vspace{12pt} \
whose legs $A,B,C$, and $D$ are attached to boundary components
$a,b,c$, and $d$ of $\S_D$ resp.
Then 
\begin{eqnarray*}
{\nabla}(F) &=& {\cases{ 
   -2 & \ \ if $\S_D$ has two boundary components and $a=b\ne c=d$, \cr
  2 & \ \ if $\S_D$ has two boundary components and $a=c\ne b=d$, \cr
  0 & \ \ otherwise.}}
 \end{eqnarray*}
\end{lemma}

\begproof
As in the proof of the previous lemma we see
that $\nabla(F)=0$ if $\S_D$ has boundary components other than $a$,
$b$, $c$, and $d$.

If there exists a component (say, $d$) to which only one leg of $K$ is
attached, then $\nabla(F)=\nabla(F_1)-\nabla(F_2)$, where
\begin{displaymath}
F_1 = 
\Picture{\DottedCircle\Arc[5]\Arc[7]\Chord[5,7]\Arc[1]\Arc[11]\Chord[1,11]
\qbezier(-0.766,-0.643)(-0.25,0)(0.66,0)
 \put(-1.25,0.68){\text{\scriptsize d}} 
 \put(-1.25,-0.77){\text{\scriptsize a}} 
 \put(0.70,0.68){\text{\scriptsize c}} 
 \put(0.75,-0.77){\text{\scriptsize b}} 
    } 
\qquad\hbox{and}\qquad
F_2 = 
\Picture{\DottedCircle\Arc[5]\Arc[7]\Arc[1]\Arc[11]\Chord[1,11]
\qbezier(-0.906,-0.423)(-0.2,0)(0.66,0)
\qbezier(-0.819,-0.574)(-0.2,0)(-0.7766,0.643)
 \put(-1.25,0.68){\text{\scriptsize d}} 
 \put(-1.25,-0.77){\text{\scriptsize a}} 
 \put(0.70,0.68){\text{\scriptsize c}} 
 \put(0.75,-0.77){\text{\scriptsize b}} 
    } 
\end{displaymath}
\\[10pt]

\noindent
and the fusion argument as above gives $\nabla(F_1)=\nabla(F_2)$
and $\nabla(F)=0$.

If $a=d$ (or $b=c$),  then $\nabla(F)=0$, since by the previous lemma
\begin{displaymath}
\nabla(F_1)=\nabla(F_2)=\mp\nabla(D')=0~,
\end{displaymath}
where 
$D'= 
\Picture{\DottedCircle\Arc[5]\Arc[7]\Chord[5,7]
 \put(-1.25,0.68){\text{\scriptsize d}} 
 \put(-1.25,-0.77){\text{\scriptsize a}} 
 \put(0.70,0.68){\text{\scriptsize c}} 
 \put(0.75,-0.77){\text{\scriptsize b}} 
    } \                     
$
\vspace{16pt}
is the chord diagram $D$ with an extra chord connecting two points
on the same boundary component of $\S_D$.

When $a=b\ne d=c$, again by the previous lemma, we have
\begin{displaymath}
\nabla(F_1)= \parbox{.8in}{\ \epsfbox{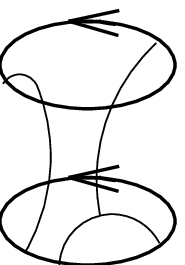}}
= \parbox{1in}{\ \epsfbox{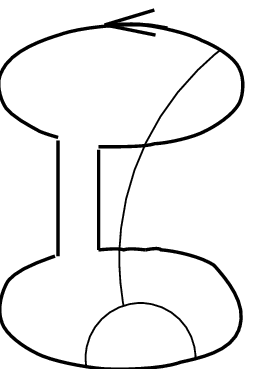}}  = -1
\end{displaymath}
and
\begin{displaymath}
\nabla(F_2)= \parbox{.8in}{\ \epsfbox{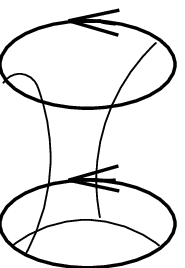}} =
\parbox{1in}{\ \epsfbox{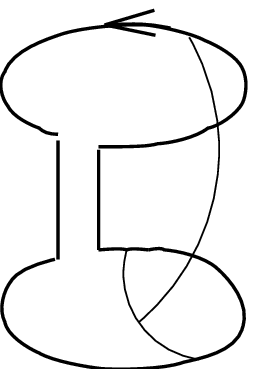}} = 1~.
\end{displaymath}
This gives $\nabla(F)=\nabla(F_1)-\nabla(F_2)=-2$.

A similar treatment of the case $a=d\ne b=c$ gives $\nabla(F)=2$.
\eproof
\bigskip

\noindent
{\em Proof of Proposition \ref{prop:rec}\/.}

Using the above lemmas we
can easily compute all five terms in both sides of
the equation (\ref{eq:conwrec}).

Denote by $p$ the number of boundary components
of the surface $\S_D$, where $D$ is the chord diagram
formed by all ``hidden'' chords in (\ref{eq:conwrec}),
and let us keep the notations of Lemma \ref{lem:2}
for the components to which the legs of the
graph $K$ are attached. 

We may  have one of the following possibilities:
\begin{enumerate}
\item If $p>2$, then  all five terms in (\ref{eq:conwrec}) vanish.
\item If  $a=d$, then by Lemma \ref{lem:4} the l.h.s. is zero, and in
the r.h.s. the second and the third terms, as well as the first and
the fourth, cancel out due to Lemma \ref{lem:2}.
\item If $a=b\ne c=d$, the last two terms of the r.h.s. vanish, and the
first two together give $2$, which is the value of the l.h.s. by
Lemma \ref{lem:4}.
\item At last, if $a=c\ne b=d$, then the l.h.s. is equal to $-2$, and
in the r.h.s. the first two terms vanish, and the last two give $-2$.
\end{enumerate}
\eproof



\appendix
\section{Values of $W_{gl(1|1)}$ for chord diagrams of
orders 4 and 5} 

In this appendix we tabulate the values of the universal weight system
$W_{gl(1|1)}$ for chord diagrams of orders 4 and 5.  These values have
been obtained by computer, after implementing the recursion relation
(\ref{eq:rec}) in {\em Mathematica}\footnote{In fact, we have computed
the values of $W_{gl(1|1)}$ on all chord diagrams of order $\leq 7$, and
they are available upon request.}.  Our results agree with the direct
computation of the weight system using the Lie superalgebraic
definition, which was implemented on a computer by Craig Snydal as
part of his undergraduate thesis \cite{Craig}.  Because of 1.3.5, only
indecomposable diagrams are tabulated.

\begin{table}[h!]
\begin{center}
\begin{tabular}{||c|l||c|l||} \hline
$D$ & \multicolumn{1}{c||}{$W_{gl(1|1)}(D)$} & $D$ &
\multicolumn{1}{c||}{$W_{gl(1|1)}(D)$}\\
\hline\hline
&&&\\ [3pt]
\Picture{
\FullCircle\EndChord[0,3]\EndChord[2,5]\EndChord[4,7]\EndChord[6,9]}
&$c^4 - 3 c^2 y + y^2$&
\Picture{
\FullCircle\EndChord[8,10]\EndChord[3,9]\EndChord[0,4]\EndChord[2,6]}
&$c^4 - 4 c^2 y + y^2$\\ [20pt]
\Picture{
\FullCircle\EndChord[0,6]\EndChord[3,9]\EndChord[1,11]\EndChord[5,7]}
&$c^4 - 3 c^2 y$&
\Picture{
\FullCircle\EndChord[1,7]\EndChord[2,10]\EndChord[4,8]\EndChord[5,11]}
&$c^4 - 5 c^2 y$\\ [20pt]
\Picture{
\FullCircle\EndChord[1,5]\EndChord[2,10]\EndChord[4,8]\EndChord[7,11]}
&$c^4 - 4 c^2 y$&
\Picture{
\FullCircle\EndChord[0,6]\EndChord[3,9]\EndChord[2,8]\EndChord[4,11]}
&$c^4 - 6 c^2 y + y^2$\\ [20pt]
\hline
\end{tabular}
\end{center}
\caption{Indecomposable chord diagrams of order $4$}
\end{table}

\clearpage

\begin{table}
\begin{center}
\begin{tabular}{||c|l||c|l||} \hline
$D$ & \multicolumn{1}{c||}{$W_{gl(1|1)}(D)$} & $D$ &
\multicolumn{1}{c||}{$W_{gl(1|1)}(D)$}\\
\hline\hline
&&&\\ [3pt]
\Picture{
\FullCircle
\EndChord[1, 3]
\EndChord[2, 5]
\EndChord[4, 7]
\EndChord[6, 9]
\EndChord[8, 10]}
&$c^5 - 4 c^3 y + 3 c y^2$
&
\Picture{
\FullCircle
\EndChord[1, 3]
\EndChord[2, 5]
\EndChord[4, 8]
\EndChord[6, 9]
\EndChord[7, 10]}
&$c^5 - 5 c^3 y + 4 c y^2$
\\ [20pt] 
\Picture{
\FullCircle
\EndChord[1, 3]
\EndChord[2, 5]
\EndChord[4, 8]
\EndChord[6, 10]
\EndChord[7, 9]}
&$c^5 - 4 c^3 y + 2 c y^2$
&
\Picture{
\FullCircle
\EndChord[1, 3]
\EndChord[2, 5]
\EndChord[4, 9]
\EndChord[6, 8]
\EndChord[7, 10]}
&$c^5 - 4 c^3 y + 3 c y^2$
\\ [20pt] 
\Picture{
\FullCircle
\EndChord[1, 3]
\EndChord[2, 6]
\EndChord[4, 7]
\EndChord[5, 9]
\EndChord[8, 10]}
&$c^5 - 5 c^3 y + 3 c y^2$
&
\Picture{
\FullCircle
\EndChord[1, 3]
\EndChord[2, 6]
\EndChord[4, 8]
\EndChord[5, 9]
\EndChord[7, 10]}
&$c^5 - 6 c^3 y + 3 c y^2$
\\ [20pt] 
\Picture{
\FullCircle
\EndChord[1, 3]
\EndChord[2, 6]
\EndChord[4, 8]
\EndChord[5, 10]
\EndChord[7, 9]}
&$c^5 - 5 c^3 y + 3 c y^2$
&
\Picture{
\FullCircle
\EndChord[1, 3]
\EndChord[2, 6]
\EndChord[4, 9]
\EndChord[5, 7]
\EndChord[8, 10]}
&$c^5 - 4 c^3 y + 2 c y^2$
\\ [20pt] 
\Picture{
\FullCircle
\EndChord[1, 3]
\EndChord[2, 6]
\EndChord[4, 9]
\EndChord[5, 8]
\EndChord[7, 10]}
&$c^5 - 5 c^3 y + 2 c y^2$
&
\Picture{
\FullCircle
\EndChord[1, 3]
\EndChord[2, 6]
\EndChord[4, 10]
\EndChord[5, 8]
\EndChord[7, 9]}
&$c^5 - 4 c^3 y + 2 c y^2$
\\ [20pt] 
\Picture{
\FullCircle
\EndChord[1, 3]
\EndChord[2, 7]
\EndChord[4, 8]
\EndChord[5, 9]
\EndChord[6, 10]}
&$c^5 - 7 c^3 y + 4 c y^2$
&
\Picture{
\FullCircle
\EndChord[1, 3]
\EndChord[2, 7]
\EndChord[4, 8]
\EndChord[5, 10]
\EndChord[6, 9]}
&$c^5 - 6 c^3 y + 2 c y^2$
\\ [20pt] 
\Picture{
\FullCircle
\EndChord[1, 3]
\EndChord[2, 7]
\EndChord[4, 9]
\EndChord[5, 8]
\EndChord[6, 10]}
&$c^5 - 6 c^3 y + 2 c y^2$
&
\Picture{
\FullCircle
\EndChord[1, 3]
\EndChord[2, 7]
\EndChord[4, 9]
\EndChord[5, 10]
\EndChord[6, 8]}
&$c^5 - 5 c^3 y + 2 c y^2$
\\ [20pt] 
\Picture{
\FullCircle
\EndChord[1, 3]
\EndChord[2, 7]
\EndChord[4, 10]
\EndChord[5, 8]
\EndChord[6, 9]}
&$c^5 - 5 c^3 y + 2 c y^2$
&
\Picture{
\FullCircle
\EndChord[1, 3]
\EndChord[2, 7]
\EndChord[4, 10]
\EndChord[5, 9]
\EndChord[6, 8]}
&$c^5 - 4 c^3 y$
\\ [20pt] 
\hline
\end{tabular}
\end{center}
\caption{Indecomposable chord diagrams of order $5$}
\end{table}

\clearpage

\begin{table}
\begin{center}
\begin{tabular}{||c|l||c|l||} \hline
$D$ & \multicolumn{1}{c||}{$W_{gl(1|1)}(D)$} & $D$ &
\multicolumn{1}{c||}{$W_{gl(1|1)}(D)$}\\
\hline\hline
&&&\\ [3pt]
\Picture{
\FullCircle
\EndChord[1, 3]
\EndChord[2, 8]
\EndChord[4, 7]
\EndChord[5, 9]
\EndChord[6, 10]}
&$c^5 - 6 c^3 y + 3 c y^2$
&
\Picture{
\FullCircle
\EndChord[1, 3]
\EndChord[2, 8]
\EndChord[4, 7]
\EndChord[5, 10]
\EndChord[6, 9]}
&$c^5 - 5 c^3 y + 2 c y^2$
\\ [20pt] 
\Picture{
\FullCircle
\EndChord[1, 3]
\EndChord[2, 9]
\EndChord[4, 7]
\EndChord[5, 8]
\EndChord[6, 10]}
&$c^5 - 5 c^3 y + 4 c y^2$
&
\Picture{
\FullCircle
\EndChord[1, 3]
\EndChord[2, 9]
\EndChord[4, 7]
\EndChord[5, 10]
\EndChord[6, 8]}
&$c^5 - 4 c^3 y + 3 c y^2$
\\ [20pt] 
\Picture{
\FullCircle
\EndChord[1, 4]
\EndChord[2, 5]
\EndChord[3, 8]
\EndChord[6, 9]
\EndChord[7, 10]}
&$c^5 - 6 c^3 y + 5 c y^2$
&
\Picture{
\FullCircle
\EndChord[1, 4]
\EndChord[2, 6]
\EndChord[3, 8]
\EndChord[5, 9]
\EndChord[7, 10]}
&$c^5 - 7 c^3 y + 3 c y^2$
\\ [20pt] 
\Picture{
\FullCircle
\EndChord[1, 4]
\EndChord[2, 6]
\EndChord[3, 9]
\EndChord[5, 8]
\EndChord[7, 10]}
&$c^5 - 6 c^3 y + 4 c y^2$
&
\Picture{
\FullCircle
\EndChord[1, 4]
\EndChord[2, 7]
\EndChord[3, 8]
\EndChord[5, 9]
\EndChord[6, 10]}
&$c^5 - 8 c^3 y + 3 c y^2$
\\ [20pt] 
\Picture{
\FullCircle
\EndChord[1, 4]
\EndChord[2, 7]
\EndChord[3, 8]
\EndChord[5, 10]
\EndChord[6, 9]}
&$c^5 - 7 c^3 y$
&
\Picture{
\FullCircle
\EndChord[1, 4]
\EndChord[2, 8]
\EndChord[3, 7]
\EndChord[5, 9]
\EndChord[6, 10]}
&$c^5 - 7 c^3 y + 2 c y^2$
\\ [20pt] 
\Picture{
\FullCircle
\EndChord[1, 4]
\EndChord[2, 8]
\EndChord[3, 7]
\EndChord[5, 10]
\EndChord[6, 9]}
&$c^5 - 6 c^3 y$
&
\Picture{
\FullCircle
\EndChord[1, 4]
\EndChord[2, 9]
\EndChord[3, 6]
\EndChord[5, 8]
\EndChord[7, 10]}
&$c^5 - 5 c^3 y + 5 c y^2$
\\ [20pt] 
\Picture{
\FullCircle
\EndChord[1, 5]
\EndChord[2, 7]
\EndChord[3, 8]
\EndChord[4, 9]
\EndChord[6, 10]}
&$c^5 - 9 c^3 y + 2 c y^2$
&
\Picture{
\FullCircle
\EndChord[1, 5]
\EndChord[2, 7]
\EndChord[3, 9]
\EndChord[4, 8]
\EndChord[6, 10]}
&$c^5 - 8 c^3 y$
\\ [20pt] 
\Picture{
\FullCircle
\EndChord[1, 6]
\EndChord[2, 7]
\EndChord[3, 8]
\EndChord[4, 9]
\EndChord[5, 10]}
&$c^5 - 10 c^3 y + 5 c y^2$
&&\\ [20pt]\hline
\end{tabular}
\end{center}
\caption{Indecomposable chord diagrams of order $5$ (cont'd)}
\end{table}

\clearpage


\end{document}